\documentclass[11pt,oneside,reqno]{amsart}
\usepackage[T1]{fontenc}
\usepackage[latin9]{inputenc}
\usepackage{color}
\usepackage{float}
\usepackage{amstext}
\usepackage{amssymb}
\usepackage{graphicx}
\usepackage[a4paper]{geometry}
\usepackage{setspace}
\usepackage[authoryear]{natbib}
\usepackage[unicode=true,pdfusetitle,
 bookmarks=true,bookmarksnumbered=false,bookmarksopen=false,
 breaklinks=true,pdfborder={0 0 0},pdfborderstyle={},backref=page,colorlinks=true]
 {hyperref}

\makeatletter
\usepackage{amsfonts}
\usepackage{amstext}

\newtheorem{thm}{Theorem}
\newtheorem{lem}{Lemma}

\newtheorem{asm}{Assumption}

\theoremstyle{definition}
\newtheorem{defn}{Definition}

\usepackage[dvipsnames]{xcolor}
\definecolor{darkblue}{rgb}{0.0, 0.0, 0.55}
\definecolor{teal}{rgb}{0.0, 0.5, 0.5}
\hypersetup{linkcolor=teal, urlcolor=blue, citecolor=darkblue}

\renewcommand*{\backrefalt}[4]{%
     \ifcase #1 %
     \or
         [Cited on page #2]
     \else
         [Cited on pages #2]
     \fi}

\usepackage{subfigure}

\makeatother

\begin{document}
\title{Optimal testing in a class of nonregular models}
\author{Yuya Shimizu}
\address{Department of Economics, University of Wisconsin, Madison, 1180 Observatory Drive, Madison, WI 53706-1393, USA.}
\email{\href{mailto:yuya.shimizu@wisc.edu}{yuya.shimizu@wisc.edu}}
\author{Taisuke Otsu}
\address{Department of Economics, London School of Economics, Houghton Street, London, WC2A 2AE, UK.}
\email{\href{mailto:t.otsu@lse.ac.uk}{t.otsu@lse.ac.uk}}
\thanks{\textit{This version: \today}}
\thanks{We are grateful to Jack Porter for his valuable comments on the manuscript. We also thank the anonymous referees for their constructive comments and suggestions. We further thank Xiaohong Chen, Harold Chiang, Greg Cox, Bruce Hansen, Ulrich M{\"u}ller, Kensuke Sakamoto, Xiaoxia Shi, and Kohei Yata for helpful comments.}
\begin{abstract}
This paper studies optimal hypothesis testing for nonregular econometric models with parameter-dependent support. We consider both one-sided and two-sided hypothesis testing and develop asymptotically uniformly most powerful tests based on a limit experiment. Our two-sided test becomes asymptotically uniformly most powerful without imposing further restrictions such as unbiasedness, and can be inverted to construct a confidence set for the nonregular parameter. Simulation results illustrate desirable finite-sample properties of the proposed tests.
\end{abstract}

\maketitle

\section{Introduction}

This paper studies optimal hypothesis testing of a class of nonregular econometric models in which the boundary of the support of the observed data depends on some parameter of interest. Such nonregular models, which typically imply discontinuous likelihood functions and nonstandard convergence rates of estimators, have often been studied in the econometrics literature; see \citet{flinn1982new}, \citet{smith1985maximum}, \citet{christensen1991exact}, \citet{donald1993maximum}, \citet{hong1998maximum}, \citet{donald2002superconsistent}, \citet{hirano2003asymptotic}, and \citet{chernozhukov2004likelihood}, among others. In contrast to most existing papers that focus on point estimation, this paper is concerned with optimal (composite) \emph{hypothesis testing} for nonregular models instead of point estimation.

For testing a simple null hypothesis against a simple alternative one, Neyman-Pearson's fundamental lemma yields an optimal power property of the likelihood ratio test even for the case of parameter-dependent support. However, the optimality result is no longer available for general testing problems with composite hypotheses. On the other hand, for regular statistical models, it is known that standard testing methods (such as the likelihood ratio, Wald, and score tests) can achieve certain asymptotic optimal power properties for testing general composite hypotheses (see Chapter 15 of \citealp{lehmann2022testing}). An open question is whether we can establish an analogous asymptotic optimality result for testing composite hypotheses on nonregular parameters in the case of parameter-dependent support, and this paper addresses this question in a positive way.

In this paper, we consider one-sided and two-sided hypothesis testing for parametric models with parameter-dependent support and develop an asymptotically uniformly most powerful (AUMP) test based on a limit experiment. Interestingly, our two-sided test attains the AUMP property without imposing further restrictions such as unbiasedness, and can be inverted to construct a confidence set for the nonregular parameter. To clarify the empirical setting, we first present a general framework and an economic application in Section \ref{sec:framework}. Then we present the main results under a benchmark setup in Section \ref{sec:bench}, where there is no covariate or nuisance parameter. Finally, we extend our optimality results to a general setup in Section \ref{sec:gen} that involves covariates and nuisance parameters. Our simulation results in Section \ref{sec:sim} illustrate desirable finite-sample properties of the proposed tests.

The most closely related papers to ours are \citet{hirano2003asymptotic} and \citet{chernozhukov2004likelihood}. \citet{hirano2003asymptotic} studied efficient point estimation of parameter-dependent support models by extending the limit of experiments argument, and showed that the Bayes estimator is asymptotically efficient under a minimax criterion but the maximum likelihood estimator is generally inefficient. To the best of our knowledge, this paper is the first one that studies optimal hypothesis testing for nonregular models with parameter-dependent support. In contrast to the Bayes estimator in \citet{hirano2003asymptotic} that involves priors on parameters, our optimal testing methods are developed based on the limiting likelihood ratio process without priors. \citet{chernozhukov2004likelihood} also investigated nonstandard asymptotic properties of estimation and testing methods for parameter-dependent support models. They established asymptotic optimality of the Bayes estimators in terms of the asymptotic average risk, and showed that the Wald test and Bayes posterior quantiles are valid for inference. However, they did not discuss optimality of the testing methods. In addition to these papers, \citeauthor{chen2018monte} (\citeyear{chen2018monte}, Appendix C) examined partially identified models with parameter-dependent support, and \citet{chen2026identification} adapted their quasi-Bayesian likelihood ratio statistics to partially identified auction models. In our simulation study, we demonstrate that the proposed test exhibits more reasonable performance than the Wald-type test of \citet{chernozhukov2004likelihood}.

\section{General framework and application}
\label{sec:framework}

We begin with a general model in which the support boundary of the conditional distribution depends on the parameter of interest:
\begin{equation}
    Y|X=x\sim f(y|x,\theta,\gamma)\mathbb{I}\{y\geq g(x,\theta)\},\label{eq:y|x,framework}
\end{equation}
where $f(\cdot|x,\theta,\gamma)$ is a smooth density on the support, $\theta$ is the scalar nonregular parameter of interest, $\gamma$ is a nuisance parameter, and $X$ is a vector of covariates. The key nonregular feature is that the boundary $g(x,\theta)$ depends on $\theta$. Thus, local changes in $\theta$ shift the support boundary.

Models with parameter-dependent support arise naturally in auction settings. As an illustration, consider the symmetric first-price procurement auction example in \citet{hirano2003asymptotic} (pp. 1320-1321). Let $X$ denote observed auction characteristics, and suppose that each bidder's cost $C$ follows a conditional exponential distribution:
\begin{equation*}
    f_{C}(c|x,\theta)=\frac{1}{s(x,\theta)}\exp\left(-\frac{c}{s(x,\theta)}\right)\mathbb{I}\{c\geq0\},
\end{equation*}
where $s(x,\theta)=E(C|x,\theta)$ is the conditional mean cost. With $m$ bidders, symmetric Bayesian-Nash equilibrium bidding implies $B=\beta(C|x,\theta)=C+s(x,\theta)/(m-1)$.
Hence the bid distribution is a shift of the exponential cost distribution:
\begin{equation*}
    f_{B}(b|x,\theta)=\frac{1}{s(x,\theta)}\exp\left(-\frac{b-\frac{s(x,\theta)}{m-1}}{s(x,\theta)}\right)\mathbb{I}\left\{ b\geq\frac{s(x,\theta)}{m-1}\right\}.
\end{equation*}
The researcher observes an iid sample of $n$ observations on $(B,X)$ and wishes to make inference on $\theta$, which characterizes the unobserved cost distribution. This model is a special case of (\ref{eq:y|x,framework}) obtained by taking $Y=B$, omitting the nuisance parameter, and setting
\begin{equation*}
    f(y|x,\theta)=\frac{1}{s(x,\theta)}\exp\left(-\frac{y-\frac{s(x,\theta)}{m-1}}{s(x,\theta)}\right),\quad \text{ and } \quad g(x,\theta)=\frac{s(x,\theta)}{m-1}.
\end{equation*}
Consequently, testing hypotheses on $\theta$ is a nonregular inference problem because $\theta$ determines the support boundary of the observed bids. The next section introduces a benchmark model to illustrate the main ideas before returning to the general framework.

\section{Benchmark case \label{sec:bench}}

This section studies a simple nonregular model with a scalar parameter and no covariates. The model provides a benchmark for the general theory developed later and includes the uniform distribution as a canonical example.

Let $Y\in\mathcal{Y}\subset\mathbb{R}$ follow the parametric model
\begin{equation}
	Y\sim f(y|\theta)\mathbb{I}\{y\geq g(\theta)\},\label{eq:f}
\end{equation}
where $\mathbb{I}\{\cdot\}$ is the indicator function, $\theta\in\Theta\subset\mathbb{R}$ is the true value of a scalar parameter, and conditions on the functions $f$ and $g$ are specified below. 
For given $\theta_{0}\in\Theta$, we consider the reparametrization $h=n(\theta-\theta_{0})\in\mathcal{H}$.
We focus on the case of $\nabla_{\theta}g(\theta_{0})>0$. The case of $\nabla_{\theta}g(\theta_{0})<0$ is analyzed in the same manner.
We consider the one-sided testing problem
\begin{equation}
    H_{0}^{-}:h\leq0\quad\text{against}\quad H_{1}^{-}:h>0\text{ (or }H_{0}^{+}:h\geq0\quad\text{against}\quad H_{1}^{+}:h<0\text{)}.\label{eq:hypo-h}
\end{equation}
We wish to test $H_{0}^{-}$ or $H_{0}^{+}$ based on an independent and identically distributed (iid) sample $Y^{n}=(Y_{1},\ldots,Y_{n})$ of $Y$.

Our testing procedure is constructed based on the limiting likelihood ratio process. Since the joint density of $Y^{n}$ is written as
\begin{equation*}
	\mathrm{d}P_{\theta}^{n}(y^{n})=\mathbb{I}\{y_{(1)}\geq g(\theta)\}\prod_{i=1}^{n}f(y_{i}|\theta),
\end{equation*}
where $y_{(1)}=\min\{y_{1},\ldots,y_{n}\}$, the likelihood ratio process on a parameter space $\mathcal{H}\subset\mathbb{R}$ is
\begin{equation*}
	Z_{n}(h,\bar{h}):=\frac{\mathrm{d}P_{\theta_{0}+\frac{\bar{h}}{n}}^{n}}{\mathrm{d}P_{\theta_{0}+\frac{h}{n}}^{n}}(Y^{n})=\frac{\mathbb{I}\{Y_{(1)}\geq g(\theta_{0}+\bar{h}/n)\}}{\mathbb{I}\{Y_{(1)}\geq g(\theta_{0}+h/n)\}}\prod_{i=1}^{n}\frac{f(Y_{i}|\theta_{0}+\bar{h}/n)}{f(Y_{i}|\theta_{0}+h/n)},
\end{equation*}
for $h,\bar{h}\in\mathcal{H}$, where the last equality holds almost surely under $P_{\theta_{0}+h/n}$.\footnote{
    We define the likelihood ratio as the density of the absolutely continuous part of
    $P_{\theta_{0}+\bar{h}/n}^{n}$ with respect to $P_{\theta_{0}+h/n}^{n}$, following \citet{van1998asymptotic}, pp.~85--86.
} To characterize asymptotic properties of this process, we impose the following assumptions.

\begin{asm}\label{asm:bench}$\quad$
    \begin{description}
        \item [{(i)}] $\{Y_{i}\}_{i=1}^{n}$ is an iid sample of $Y\in\mathcal{Y}\subset\mathbb{R}$ with the Lebesgue density in (\ref{eq:f}). The parameter space $\Theta$ is convex, $\theta_{0}\in\mathrm{int}(\Theta)$, and $0\in\mathcal{H}$.
        \item [{(ii)}] $f(y|\theta)$ is twice continuously differentiable in $\theta$ for all $y$. In some open neighborhood $\mathcal{N}$ of $\theta_{0}$, $f(y|\theta)$ and $\nabla_{\theta}f(y|\theta)$ are continuous in $y$ for $\theta\in\mathcal{N}$, there exists a constant $C$ such that $0<f(y|\theta)<C<\infty$ and $|\nabla_{\theta}f(y|\theta)|<C<\infty$ for all $y$ and $\theta\in\mathcal{N}$, and
			\begin{align*}
				\int\sup_{\theta\in\mathcal{N}}\left\Vert \nabla_{\theta}f(y|\theta)\right\Vert \mathbb{I}\{y\geq g(\theta)\}dy &<\infty,\\
				\int\sup_{\theta,\bar{\theta}\in\mathcal{N}}\frac{\left\Vert \nabla_{\theta}f(y|\bar{\theta})\right\Vert ^{2}}{f(y|\bar{\theta})^{2}}\mathbb{I}\{y\geq g(\theta)\}f(y|\theta)dy &<\infty,\\
				\int\sup_{\theta,\bar{\theta}\in\mathcal{N}}\frac{\left\Vert \nabla_{\theta\theta}f(y|\bar{\theta})\right\Vert }{f(y|\bar{\theta})}\mathbb{I}\{y\geq g(\theta)\}f(y|\theta)dy &<\infty.
			\end{align*}
              $g(\theta)$ is continuously differentiable in $\theta$, $\nabla_{\theta}g(\theta_{0})>0$, and $\sup_{\theta\in\mathcal{N}}\left\Vert \nabla_{\theta}g(\theta)\right\Vert <\infty$.
    \end{description}
\end{asm}

Assumption \ref{asm:bench} (i) is standard, and Assumption \ref{asm:bench} (ii) contains smoothness and boundedness conditions on the functions $f$ and $g$ in (\ref{eq:f}). Assumption \ref{asm:bench} guarantees weak convergence (denoted by ``$\rightsquigarrow$'') of the likelihood ratio process.
\begin{lem}\label{lem:weak_conv}
	Under Assumption \ref{asm:bench}, it holds that
    \begin{equation}
    \{Z_{n}(h,\bar{h})\}_{\bar{h}\in I}\rightsquigarrow\{Z(h,\bar{h})\}_{\bar{h}\in I}:=\{e^{(\bar{h}-h)/\lambda}D_{h,\bar{h}}\}_{\bar{h}\in I}\quad\text{under }P_{\theta_{0}+\frac{h}{n}},\label{eq:Z}
    \end{equation}
    for every finite $I\subset\mathcal{H}$, where $\lambda=\{f(g(\theta_{0})|\theta_{0})\nabla_{\theta}g(\theta_{0})\}^{-1}$ and
    \begin{equation}
        D_{h,\bar{h}}:=\mathbb{I}\{W_{h}>\bar{h}\}\text{ with }W_{h}\sim f_{W}(w|h)=\frac{1}{\lambda}e^{-(w-h)/\lambda}\mathbb{I}\{w>h\}.
        \label{eq:W}
    \end{equation}
    Moreover, we can show the convergences for the components as
	\begin{equation}
		\mathbb{I}\left\{Y_{(1)}\geq g\left(\theta_{0}+\frac{\bar{h}}{n}\right)\right\}
		\rightsquigarrow D_{h,\bar{h}}
		\quad\text{under }P_{\theta_{0}+\frac{h}{n}},
		\label{eq:D}
	\end{equation}
	and
	\begin{equation}
		\prod_{i=1}^{n}
		\frac{f(Y_{i}|\theta_{0}+\bar{h}/n)}
		{f(Y_{i}|\theta_{0}+h/n)}
		\stackrel{p}{\rightarrow}
		e^{(\bar{h}-h)/\lambda}
		\quad\text{under }P_{\theta_{0}+\frac{h}{n}},
		\label{eq:exp}
	\end{equation}
	for every $\bar{h}\in\mathcal{H}$.
\end{lem}
The limiting likelihood ratio process has several notable features. Note that $\lambda$ is a known constant in the present setup. In contrast to likelihood ratio processes arising in locally asymptotically normal models, $Z_{n}(h,\bar{h})$ converges to a limit whose behavior is determined by the binary variable $D_{h,\bar{h}}$. This is due to the lack of differentiability in quadratic mean of the density $\mathrm{d}P_{\theta}^{n}(y^{n})$. If $\bar{h}>h$, then $D_{h,\bar{h}}$ is nondegenerate, and the distribution function of $Z(h,\bar{h})$ has jump discontinuities at $z=0$ and $z=e^{(\bar{h}-h)/\lambda}$. If $\bar{h}\leq h$, then $D_{h,\bar{h}}=1$ almost surely, and the distribution function has a single jump at $z=e^{(\bar{h}-h)/\lambda}$. Thus, we cannot use the conventional asymptotic theory based on a quadratic expansion of the likelihood ratio process to evaluate asymptotic size and power.

Even though the likelihood ratio process $\{Z_{n}(h,\bar{h})\}_{\bar{h}\in I}$ exhibits such nonregularity, the family of distributions in the limit experiment defined in (\ref{eq:W}) has monotone likelihood ratio in $W$.

\begin{defn}[Monotone likelihood ratio]\label{def:MLR}
	(\citealp{shao2003mathematical}, Definition 6.2) Suppose that the distribution of $W$ is in $\mathcal{P}=\{P_{h}:h\in\mathcal{H}\}$, a parametric family indexed by a real-valued $h$, and that $\mathcal{P}$ is dominated by a $\sigma$-finite measure $\mu$. Let
	$p_{h}=\mathrm{d}P_{h}/\mathrm{d}\mu$.
    The family $\mathcal{P}$ is said to have \textit{monotone likelihood ratio} in the real-valued statistic $S(W)$ if, for every $h_{1}<h_{2}$, the ratio $p_{h_{2}}(w)/p_{h_{1}}(w)$ is a nondecreasing function of $S(w)$ on the set where at least one of $p_{h_{1}}(w)$ and $p_{h_{2}}(w)$ is positive.
\end{defn}

A key feature of distributions with monotone likelihood ratio is the existence of the uniformly most powerful (UMP) test. Relying on the limiting likelihood ratio process is essential since the finite-sample likelihood ratio does not exhibit a monotone likelihood ratio property in general.

\begin{lem}\label{lem:MLRT} (\citealp{shao2003mathematical}, Theorem 6.2.1) Suppose that a random variable $U_{h}$ has a distribution in $\mathcal{P}=\{P_{h}:h\in\mathcal{H}\subset\mathbb{R}\}$ that has monotone likelihood ratio in $S(U_{h})$. Consider the problem of testing $H_{0}:h\leq h_{0}$ against $H_{1}:h>h_{0}$, where $h_{0}$ is a given constant. Then there exists a UMP test of size $\alpha$ given by
    \begin{equation}
        \phi(U_{h})=\begin{cases}
            1      & \text{ if }S(U_{h})>c \\
            \kappa & \text{ if }S(U_{h})=c \\
            0      & \text{ if }S(U_{h})<c
        \end{cases},\label{eq:UMP}
    \end{equation}
    where $c$ and $\kappa$ are determined by $E_{h_{0}}[\phi(U_{h_{0}})]=\alpha$. \end{lem}

Here $\phi(U_{h})=1$ and $0$ mean rejection and acceptance of $H_{0}$, respectively, and $\phi(U_{h})=\kappa$ means rejection with probability $\kappa$. Therefore, to derive an AUMP test for $H_{0}^{-}:h\leq0$, we can still invoke the asymptotic representation lemma below to argue that the sample counterpart of the monotone likelihood ratio test with $S(W_{h})=W_{h}$ is AUMP.

\begin{lem}\label{lem:vdv} (\citealp{van1998asymptotic}, Theorem 15.1) Let the sequence of experiments $\mathcal{E}_{n}=\{P_{n,h}:h\in\mathcal{H}\}$ converge to a dominated experiment $\mathcal{E}=\{P_{h}:h\in\mathcal{H}\}$. Suppose that a sequence of power functions $\pi_{n}$ of tests in $\mathcal{E}_{n}$ converges pointwise, i.e., $\pi_{n}(h)\rightarrow\pi(h)$ for every $h$ and some function $\pi$. Then $\pi$ is a power function in the limit experiment, i.e., there exists a test $\phi$ in $\mathcal{E}$ with $\pi(h)=E_{h}[\phi(X)]$ for every $h$. \end{lem}

To derive an AUMP test for the other one-sided hypothesis $H_{0}^{+}:h\geq0$, we utilize the fact that the limit experiment can be represented as a function of the random variable $W_{h}$. Accordingly, we apply the Neyman-Pearson lemma (e.g., Theorem 3.2.1 (ii) of \citealp{lehmann2022testing}) to $W_{h}$.

Hereafter, we first formalize this argument for one-sided tests on $H_{0}^{-}:h\leq0$ and $H_{0}^{+}:h\geq0$ (Section \ref{sub:os}), and then extend the argument to two-sided testing for $H_{0}:h=0$ (Section \ref{sub:ts}).

\subsection{One-sided tests\label{sub:os}}

First, we consider testing $H_{0}^{-}:h\leq0$ against $H_{1}^{-}:h>0$ at the significance level $\alpha\in(0,1)$. By taking a sample counterpart of $\phi(W_{h})$ in (\ref{eq:UMP}) with the limiting process $W_{h}$ in (\ref{eq:W}), our test is constructed as
\begin{equation}
    \phi_{n}^{-}(Y^{n})=\left\{ \begin{array}{lll}
        1 & \text{if} & Y_{(1)}>g\left(\theta_{0}+\frac{\lambda}{n}\log\left(\frac{1}{\alpha}\right)\right)     \\
        0 & \text{if} & Y_{(1)}\leq g\left(\theta_{0}+\frac{\lambda}{n}\log\left(\frac{1}{\alpha}\right)\right)
    \end{array}\right..\label{eq:phi-}
\end{equation}
This test achieves the asymptotic optimality property in Definition \ref{def:AUMP} below, introduced by \citet{choi1996asymptotically}.

\begin{defn}[Asymptotically uniformly most powerful test]\label{def:AUMP} For testing $H_{0}:h\leq0$ against $H_{1}:h>0$ (or $H_{0}:h\geq0$ against $H_{1}:h<0$ or $H_{0}:h=0$ against $H_{1}:h\neq0$), a sequence of tests $\{\phi_{n}\}$ is called \textit{asymptotically uniformly most powerful (AUMP)} in $\mathcal{H}$ at asymptotic level $\alpha$ if $\limsup_{n}E_{\theta_{0}+h/n}[\phi_{n}]\leq\alpha$ for every $h\leq0$ (or $h\geq0$ or $h=0$) in $\mathcal{H}$ and for every other sequence of test functions $\{\psi_{n}\}$ satisfying $\limsup_{n}E_{\theta_{0}+h/n}[\psi_{n}]\leq\alpha$ for every $h\leq0$ (or $h\geq0$ or $h=0$) in $\mathcal{H}$,
	\begin{equation*}
		\liminf_{n}E_{\theta_{0}+h/n}[\phi_{n}]\geq\limsup_{n}E_{\theta_{0}+h/n}[\psi_{n}],
	\end{equation*}
    for every $h>0$ (or $h<0$ or $h\neq0$) in $\mathcal{H}$.\end{defn}

The following asymptotic optimality result is established using the monotone likelihood ratio property of the limit experiment, as demonstrated in the Appendix.

\begin{thm}\label{thm:opt-} Suppose that Assumption \ref{asm:bench} holds for every $h\in\mathcal{H}$. Then the test $\phi_{n}^{-}(Y^{n})$ is AUMP in $\mathcal{H}$ at level $\alpha$ for testing $H_{0}^{-}:h\leq0$ against $H_{1}^{-}:h>0$. \end{thm}

Next, we consider the other one-sided testing problem $H_{0}^{+}:h\geq0$ against $H_{1}^{+}:h<0$. The basic idea is the same as in the previous case. We propose the following test:
\begin{equation}
    \phi_{n}^{+}(Y^{n})=\left\{ \begin{array}{lll}
        1 & \text{if} & Y_{(1)}<\max\left\{ g(\theta_{0}),g\left(\theta_{0}+\frac{\lambda}{n}\log\left(\frac{1}{1-\alpha}\right)\right)\right\}    \\
        0 & \text{if} & Y_{(1)}\geq\max\left\{ g(\theta_{0}),g\left(\theta_{0}+\frac{\lambda}{n}\log\left(\frac{1}{1-\alpha}\right)\right)\right\}
    \end{array}\right.,\label{eq:phi+}
\end{equation}
The asymptotic optimality of this test is established by applying the Neyman-Pearson lemma to $W_{h}$.

\begin{thm}\label{thm:opt+} Suppose that Assumption \ref{asm:bench} holds for every $h\in\mathcal{H}$. Then the test $\phi_{n}^{+}(Y^{n})$ is AUMP in $\mathcal{H}$ at level $\alpha$ for testing $H_{0}^{+}:h\geq0$ against $H_{1}^{+}:h<0$. \end{thm}

Since we assume $\nabla_{\theta}g(\theta_{0})>0$ and $\lambda>0$, we have $g(\theta_{0})<g\left(\theta_{0}+\frac{\lambda}{n}\log\left(\frac{1}{1-\alpha}\right)\right)$ eventually. However, this inequality can be violated in finite samples. We can show that every test $\psi_{n}^{+}$ rejecting the null with probability one if $Y_{(1)}<g(\theta_{0})$ and rejecting the null with probability $\alpha$ if $Y_{(1)}\geq g(\theta_{0})$ is AUMP at level $\alpha$ for testing $H_{0}^{+}:h\geq0$ against $H_{1}^{+}:h<0$. We recommend using (\ref{eq:phi+}) since it avoids randomization.

\subsection{Two-sided test\label{sub:ts}}

This subsection considers two-sided testing $H_{0}:h=0$ against $H_{1}:h\neq0$. By combining the optimal one-sided tests derived in the last subsection, we propose the following (unequal-tailed) two-sided test:
\begin{equation}
    \phi_{n}(Y^{n})=\left\{ \begin{array}{lll}
        1 & \text{if} & Y_{(1)}>g\left(\theta_{0}+\frac{\lambda}{n}\log\left(\frac{1}{\alpha}\right)\right)\text{ or }Y_{(1)}<g(\theta_{0}) \\
        0 & \text{if} & g(\theta_{0})\leq Y_{(1)}\leq g\left(\theta_{0}+\frac{\lambda}{n}\log\left(\frac{1}{\alpha}\right)\right)
    \end{array}\right..\label{eq:phi_ts}
\end{equation}
Indeed, this test is shown to be AUMP.

\begin{thm}\label{thm:opt_ts} Suppose that Assumption \ref{asm:bench} holds for every $h\in\mathcal{H}$. Then the proposed test $\phi_{n}(Y^{n})$ is AUMP in $\mathcal{H}$ at level $\alpha$ for testing $H_{0}:h=0$ against $H_{1}:h\neq0$. \end{thm}

Importantly, the proposed test achieves the AUMP property without imposing additional restrictions such as unbiasedness. This feature is shared by a finite-sample UMP two-sided test for the uniform distribution (e.g., \citeauthor{lehmann2022testing}, \citeyear{lehmann2022testing}, Problem 3.2, p.~105).

Since the two-sided test $\phi_{n}(Y^{n})$ does not randomize, we can easily construct a $100(1-\alpha)$\% confidence set by the test inversion:
\begin{equation*}
	CS=\left\{ \theta\in\Theta:g(\theta)\leq Y_{(1)}\leq g\left(\theta+\frac{\lambda(\theta)}{n}\log\left(\frac{1}{\alpha}\right)\right)\right\},
\end{equation*}
where $\lambda(\theta)=\{f(g(\theta)|\theta)\nabla_{\theta}g(\theta)\}^{-1}$.
This confidence set is also pointwise optimal (i.e., asymptotically uniformly most accurate; see \citealp{choi1996asymptotically}). The next section extends these benchmark results to the general model with covariates and nuisance parameters.

\section{General case \label{sec:gen}}

In this section, as in \citet{hirano2003asymptotic}, we analyze the general framework in Section \ref{sec:framework} that accommodates discrete covariates and nuisance parameters:
\begin{equation}
	Y|X=x\sim f(y|x,\theta,\gamma)\mathbb{I}\{y\geq g(x,\theta)\},\label{eq:y|x,ga}
\end{equation}
where $\theta\in\Theta\subset\mathbb{R}$ is a scalar parameter of interest, $\gamma\in\Gamma\subset\mathbb{R}^{d}$ is a $d$-dimensional vector of (regular) nuisance parameters, $Y$ is a scalar dependent variable, and $X$ is a $d_{X}$-dimensional vector of discrete covariates with support $\mathcal{X}=\{a_{1},\dots,a_{L}\}$. On reparametrization $h=n(\theta-\theta_{0})$, this section considers the one-sided testing problem
\begin{equation*}
    H_{0}^{-}:h\leq0\quad\text{against}\quad H_{1}^{-}:h>0\text{ (or }H_{0}^{+}:h\geq0\quad\text{against}\quad H_{1}^{+}:h<0\text{)},
\end{equation*}
with the asymptotic level of significance $\alpha$ and a given $\gamma\in\Gamma$.

Let $(Y^{n},X^{n})=((Y_{1},\ldots,Y_{n}),(X_{1},\ldots,X_{n}))$ be an iid sample of $(Y,X)\in\mathbb{R}\times\mathcal{X}$. To extend our benchmark results in the last section, we consider the likelihood ratio process
\begin{equation}
	Z_{n}(h,\bar{h},\gamma):=\prod_{i=1}^{n}\frac{f(Y_{i}|X_{i},\theta_{0}+\bar{h}/n,\gamma)\mathbb{I}\{Y_{i}\geq g(X_{i},\theta_{0}+\bar{h}/n)\}}{f(Y_{i}|X_{i},\theta_{0}+h/n,\gamma)\mathbb{I}\{Y_{i}\geq g(X_{i},\theta_{0}+h/n)\}},\label{eq:LR1}
\end{equation}
for $h,\bar{h}\in\mathcal{H}$ with a parameter space $\mathcal{H}\subset\mathbb{R}$. The ratio is almost surely well-defined under the true data-generating process $P_{\theta_{0}+h/n,\gamma}$. We also consider the plug-in likelihood ratio process $Z_{n}(h,\hat{h}_{n},\hat{\gamma}_{n})$, where $\hat{\gamma}_{n}$ and $\hat{h}_{n}$ are estimators of $\gamma$ and $\bar{h}$, respectively. In contrast to the benchmark case, where the threshold for testing takes the form $g(\theta_{0}+\bar{h}/n)$, the optimal values of $\bar{h}$ for our test depend on certain population quantities that must be estimated (see Theorems \ref{thm:opt-nui} and \ref{thm:opt+nui} below). Thus, we introduce an estimator $\hat{h}_{n}$ for $\bar{h}$ in this general case.

In this section, we impose the following assumptions.

\begin{asm}\label{asm:nui}$\quad$
	\begin{description}
		\item [{(i)}] $\{(Y_{i},X_{i})\}_{i=1}^{n}$ is an iid sample taking values in $\mathcal{Y}\times\mathcal{X}$, where $\mathcal{Y}\subset\mathbb{R}$. The covariate support is $\mathcal{X}=\{a_{1},\ldots,a_{L}\}$. The conditional density of $Y$ given $X$ is specified in (\ref{eq:y|x,ga}). The parameter space $\Theta\times\Gamma$ of $(\theta,\gamma)$ is convex, $\theta_{0}\in\mathrm{int}(\Theta)$, $\gamma\in\mathrm{int}(\Gamma)$, and $0\in\mathcal{H}$.
		Moreover, the marginal distribution of $X$ satisfies $\Pr\{X=a_{j}\}>0$ for each $j=1,\ldots,L$ and does not depend on $(\theta,\gamma)$.
		\item [{(ii)}] $\{\hat{h}_{n}\}$ is a random sequence satisfying $\sqrt{n}(\hat{h}_{n}-\bar{h})=O_{P_{\theta_{0}+\frac{h}{n},\gamma}}(1)$ for some user-specified constant $\bar{h}$. $\{\hat{\gamma}_{n}\}$ is a random sequence satisfying $\sqrt{n}(\hat{\gamma}_{n}-\gamma)=O_{P_{\theta_{0}+\frac{h}{n},\gamma}}(1)$.
		\item [{(iii)}] Let $\beta=(\theta,\gamma^{\prime})^{\prime}$ and $\beta_{0}=(\theta_{0},\gamma^{\prime})^{\prime}$. $f(y|x,\beta)$ is twice continuously differentiable in $\beta$ for all $y$ and $x$, and $g(x,\theta)$ is continuously differentiable in $\theta$ for all $x$. In some open neighborhood $\mathcal{N}$ of $\beta_{0}$, $f(y|x,\beta)$ and $\nabla_{\beta}f(y|x,\beta)$ are continuous in $y$ for $\beta\in\mathcal{N}$, there exists a constant $C$ such that $0<f(y|x,\beta)<C<\infty$ and $\|\nabla_{\beta}f(y|x,\beta)\|<C<\infty$ for all $y$, $x$, and $\beta\in\mathcal{N}$, and for each $j=1,\ldots,L$,
		\begin{align*}
			&\int\sup_{\beta\in\mathcal{N}}\left\Vert\nabla_{\beta}f(y|a_{j},\beta)\right\Vert\mathbb{I}\{y\geq g(a_{j},\theta)\}dy<\infty,\\
			&\int\sup_{\bar{\beta},\beta\in\mathcal{N}}\frac{\left\Vert\nabla_{\beta}f(y|a_{j},\bar{\beta})\right\Vert^{2}}{f(y|a_{j},\bar{\beta})^{2}}\mathbb{I}\{y\geq g(a_{j},\theta)\}f(y|a_{j},\beta)dy<\infty,\\
			&\int\sup_{\bar{\beta},\beta\in\mathcal{N}}\frac{\left\Vert\nabla_{\beta\beta}f(y|a_{j},\bar{\beta})\right\Vert}{f(y|a_{j},\bar{\beta})}\mathbb{I}\{y\geq g(a_{j},\theta)\}f(y|a_{j},\beta)dy<\infty.
		\end{align*}
		$\nabla_{\theta}g(a_{j},\theta_{0})>0$ and $\sup_{\beta\in\mathcal{N}}\left\Vert\nabla_{\theta}g(a_{j},\theta)\right\Vert<\infty$ for each $j=1,\ldots,L$.
	\end{description}
\end{asm}

As specified in Assumption \ref{asm:nui} (i), we focus on the case where covariates are discrete variables as in \citet{hirano2003asymptotic}. Assumption \ref{asm:nui} (ii) imposes standard $\sqrt{n}$-rate conditions for the regular nuisance parameter and the alternative value. For example, $\gamma$ may be estimated by maximum likelihood, while $\hat{h}_{n}$ may be constructed by plugging estimators of the relevant population quantities into the expression for $\bar{h}$. Assumption \ref{asm:nui} (iii) lists boundedness and smoothness conditions on the functions $f$ and $g$.

To present the limiting distribution of the plug-in likelihood ratio process $\{Z_{n}(h,\hat{h}_{n},\hat{\gamma}_{n})\}$, we introduce further notations. Define
\begin{equation*}
	\lambda=\left\{ E_{X}[f(g(X,\theta_{0})|X,\beta_{0})\nabla_{\theta}g(X,\theta_{0})]\right\}^{-1},
\end{equation*}
$G_{j}=\nabla_{\theta}g(a_{j},\theta_{0})$, and $\lambda_{j}=\{\Pr(X=a_{j})f(g(a_{j},\theta_{0})|a_{j},\beta_{0})\}^{-1}$. Also, let $(W_{h,1},\ldots,W_{h,L})$ be mutually independent random variables with densities
\begin{equation*}
	f_{W_{j}}(w\mid h,\gamma)=e^{-(w-G_{j}h)/\lambda_{j}}\mathbb{I}\{w>G_{j}h\}/\lambda_{j}.
\end{equation*}
For every $h,\bar{h}\in\mathcal{H}$, define $D_{h,\bar{h}}=\prod_{j=1}^{L}\mathbb{I}\{W_{h,j}>G_{j}\bar{h}\}$, which is a random variable following a Bernoulli distribution
\begin{align*}
	D_{h,\bar{h}} \sim \mathrm{Bernoulli}(\exp(-E_{X}[\mathbb{I}\{\bar{h}-h\geq0\}f(g(X,\theta_{0})|X,\beta_{0})\nabla_{\theta}g(X,\theta_{0})(\bar{h}-h)])).
\end{align*}
The weak convergence of the plug-in likelihood ratio process is established as follows.

\begin{lem}\label{lem:weak_conv_nui}
	Under Assumption \ref{asm:nui}, it holds that
	\begin{equation}
		\{Z_{n}(h,\bar{h},\gamma)\}_{\bar{h}\in I}\rightsquigarrow\{Z(h,\bar{h})\}_{\bar{h}\in I}:=\{e^{(\bar{h}-h)/\lambda}D_{h,\bar{h}}\}_{\bar{h}\in I}\quad\text{under }P_{\theta_{0}+\frac{h}{n},\gamma},\label{eq:Z_ga}
	\end{equation}
	for every finite $I\subset\mathcal{H}$ and
    \begin{equation}
		Z_{n}(h,\hat{h}_{n},\hat{\gamma}_{n})\rightsquigarrow Z(h,\bar{h})\quad\text{under }P_{\theta_{0}+\frac{h}{n},\gamma},
        \label{eq:Z_ga_hat}
	\end{equation}
    for every $\bar{h}\in\mathcal{H}$ and $\hat{h}_{n}$ satisfying $\sqrt{n}(\hat{h}_{n}-\bar{h})=O_{P_{\theta_{0}+h/n,\gamma}}(1)$, where $\lambda=(\sum_{j=1}^{L}G_{j}/\lambda_{j})^{-1}$ and $D_{h,\bar{h}}=\prod_{j=1}^{L}\mathbb{I}\{W_{h,j}>G_{j}\bar{h}\}$. Moreover, we can show the convergences for the components as
	\begin{equation}
		\prod_{i=1}^{n}\mathbb{I}\{Y_{i}\geq g(X_{i},\theta_{0}+\hat{h}_{n}/n)\}\rightsquigarrow D_{h,\bar{h}}\quad\text{under }P_{\theta_{0}+\frac{h}{n},\gamma},\label{eq:D_ga}
	\end{equation}
	and
	\begin{equation}
		\prod_{i=1}^{n}\frac{f(Y_{i}|X_{i},\theta_{0}+\hat{h}_{n}/n,\hat{\gamma}_{n})}{f(Y_{i}|X_{i},\theta_{0}+h/n,\hat{\gamma}_{n})}\stackrel{p}{\rightarrow}e^{(\bar{h}-h)/\lambda}\quad\text{under }P_{\theta_{0}+\frac{h}{n},\gamma},\label{eq:exp_ga}
	\end{equation}
	for every $\bar{h}\in\mathcal{H}$.
\end{lem}

Despite the presence of covariates and nuisance parameters, the limiting likelihood ratio retains the same basic form as in the benchmark case: it is the product of a deterministic exponential term and a Bernoulli random variable. There are, however, three important differences. First, the plug-in process contains the estimators $\hat{h}_{n}$ and $\hat{\gamma}_{n}$. Second, the Bernoulli component is generated by the $L$-dimensional vector $(W_{h,1},\ldots,W_{h,L})$. Third, the limit experiment depends on the nuisance parameter $\gamma$ through $\lambda$ and the distributions of the $W_{h,j}$'s.

Hereafter, we separately consider testing $H_{0}^{-}:h\leq0$ against $H_{1}^{-}:h>0$ for every fixed $\gamma\in\Gamma$ and $H_{0}^{+}:h\geq0$ against $H_{1}^{+}:h<0$ for every fixed $\gamma\in\Gamma$. Let $\hat{\lambda}_{n}$ be a $\sqrt{n}$-consistent estimator of $\lambda$. For example, $\lambda$ can be estimated by $\hat{\lambda}_{n}=(n^{-1}\sum_{i=1}^{n}f(g(X_{i},\theta_{0})|X_{i},\theta_{0},\hat{\gamma}_{n})\nabla_{\theta}g(X_{i},\theta_{0}))^{-1}$.
Under Assumption \ref{asm:nui} (i) and (iii), it is straightforward to show that this $\hat{\lambda}_{n}$ satisfies $\sqrt{n}(\hat{\lambda}_{n}-\lambda)=O_{P_{\theta_{0}+h/n,\gamma}}(1)$ for every fixed $h\in\mathcal{H}$.
We will show the optimality defined as follows for the proposed tests in the presence of nuisance parameters.
\begin{defn}[Asymptotically uniformly most powerful test at $\gamma$]\label{def:AUMPnui}
	For testing $H_{0}:h\leq0$ against $H_{1}:h>0$ (or $H_{0}:h\geq0$ against $H_{1}:h<0$ or $H_{0}:h=0$ against $H_{1}:h\neq0$) for fixed $\gamma\in\Gamma$, a sequence of tests $\{\phi_{n}\}$ is called \textit{asymptotically uniformly most powerful (AUMP)} in $\mathcal{H}$ at asymptotic level $\alpha$ for the nuisance parameter value $\gamma$ if $\limsup_{n}E_{\theta_{0}+h/n,\gamma}[\phi_{n}]\leq\alpha$ for every $h\leq0$ (or $h\geq0$ or $h=0$) in $\mathcal{H}$ and for every other sequence of test functions $\{\psi_{n}\}$ satisfying $\limsup_{n}E_{\theta_{0}+h/n,\gamma}[\psi_{n}]\leq\alpha$ for every $h\leq0$ (or $h\geq0$ or $h=0$) in $\mathcal{H}$,
	\begin{equation*}
		\liminf_{n}E_{\theta_{0}+h/n,\gamma}[\phi_{n}]\geq\limsup_{n}E_{\theta_{0}+h/n,\gamma}[\psi_{n}],
	\end{equation*}
	for every $h>0$ (or $h<0$ or $h\neq0$) in $\mathcal{H}$.
\end{defn}

\subsection{One-sided tests\label{sub:os_nui}}

First, we consider one-sided testing $H_{0}^{-}:h\leq0$ against $H_{1}^{-}:h>0$ at the significance level $\alpha\in(0,1)$ for a fixed $\gamma\in\Gamma$. By extending the one-sided test in (\ref{eq:phi-}), our test is defined as
\begin{equation}
	\phi_{n}^{-}(\hat{h}_{n}^{-},Y^{n},X^{n})=\left\{\begin{array}{lll}
	1&\text{if}&Y_{i}>g(X_{i},\theta_{0}+\hat{h}_{n}^{-}/n)\text{ for all }i\\
	0&\text{if}&Y_{i}\leq g(X_{i},\theta_{0}+\hat{h}_{n}^{-}/n)\text{ for some }i
	\end{array}\right.,\label{eq:phi-nui}
\end{equation}
where $\hat{h}_{n}^{-}=\hat{\lambda}_{n}\log\left(1/\alpha\right)$. The idea for constructing this test is essentially the same as that in the benchmark case in Section \ref{sec:bench}. The main difference is that $\lambda$ is unknown and needs to be estimated by plugging in the $\sqrt{n}$-consistent estimator $\hat{\lambda}_{n}$. The next theorem shows that this test achieves an asymptotic optimality property.

\begin{thm}\label{thm:opt-nui}
		Suppose that Assumption \ref{asm:nui} holds for every $h\in\mathcal{H}$. Let $\{\hat{h}_{n}^{-}\}$ be a random sequence such that $\sqrt{n}(\hat{h}_{n}^{-}-\bar{h}^{-})=O_{P_{\theta_{0}+\frac{h}{n},\gamma}}(1)$ for every fixed $h\in\mathcal{H}$ and $\bar{h}^{-}=\lambda\log(1/\alpha)\in\mathcal{H}$. Then the test $\phi_{n}^{-}(\hat{h}_{n}^{-},Y^{n},X^{n})$ is AUMP in $\mathcal{H}$ at level $\alpha$ for the nuisance parameter value $\gamma$ when testing $H_{0}^{-}$ against $H_{1}^{-}$.
\end{thm}

Next, we consider the other one-sided testing problem $H_{0}^{+}:h\geq0$ against $H_{1}^{+}:h<0$ at the significance level $\alpha\in(0,1)$ for a fixed $\gamma\in\Gamma$. In this case, our test is defined as
\begin{equation}
	\phi_{n}^{+}(\hat{h}_{n}^{+},Y^{n},X^{n})=\left\{\begin{array}{lll}
	1&\text{if}&Y_{i}<\max\{g(X_{i},\theta_{0}),g(X_{i},\theta_{0}+\hat{h}_{n}^{+}/n)\}\text{ for some }i\\
	0&\text{if}&Y_{i}\geq\max\{g(X_{i},\theta_{0}),g(X_{i},\theta_{0}+\hat{h}_{n}^{+}/n)\}\text{ for all }i
	\end{array}\right.,\label{eq:phi+nui}
\end{equation}
where $\hat{h}_{n}^{+}=\hat{\lambda}_{n}\log\left(1/(1-\alpha)\right)$. As in Theorem \ref{thm:opt-nui}, the asymptotic optimality of this test is obtained as follows.

\begin{thm}\label{thm:opt+nui}
		Suppose that Assumption \ref{asm:nui} holds for every $h\in\mathcal{H}$. Let $\{\hat{h}_{n}^{+}\}$ be a random sequence such that $\sqrt{n}(\hat{h}_{n}^{+}-\bar{h}^{+})=O_{P_{\theta_{0}+\frac{h}{n},\gamma}}(1)$ for every fixed $h\in\mathcal{H}$ and $\bar{h}^{+}=\lambda\log\left(1/(1-\alpha)\right)\in\mathcal{H}$. Then the test $\phi_{n}^{+}(\hat{h}_{n}^{+},Y^{n},X^{n})$ is AUMP in $\mathcal{H}$ at level $\alpha$ for the nuisance parameter value $\gamma$ when testing $H_{0}^{+}$ against $H_{1}^{+}$.
\end{thm}

\subsection{Two-sided test\label{sub:ts_nui}}

This subsection considers two-sided testing $H_{0}:h=0$ against $H_{1}:h\neq0$ at the significance level $\alpha\in(0,1)$ for a fixed $\gamma\in\Gamma$. By combining the optimal one-sided tests derived in the last subsection, we propose the following (unequal-tailed) two-sided test:
\begin{equation}
	\phi_{n}(\hat{h}_{n}^{-},Y^{n},X^{n})=\left\{\begin{array}{lll}
		1&\text{if}&Y_{i}>g(X_{i},\theta_{0}+\hat{h}_{n}^{-}/n)\text{ for all }i\text{ or }Y_{i}<g(X_{i},\theta_{0})\text{ for some }i\\
		0&\text{if}&Y_{i}\leq g(X_{i},\theta_{0}+\hat{h}_{n}^{-}/n)\text{ for some }i\text{ and }Y_{i}\geq g(X_{i},\theta_{0})\text{ for all }i
	\end{array}\right..\label{eq:phi_ts-nui}
\end{equation}

\begin{thm}\label{thm:opt_ts-nui} Suppose that Assumption \ref{asm:nui} holds for every $h\in\mathcal{H}$. Let $\{\hat{h}_{n}^{-}\}$ be a random sequence defined in Theorem $\ref{thm:opt-nui}$. Then the test $\phi_{n}(\hat{h}_{n}^{-},Y^{n},X^{n})$ is AUMP in $\mathcal{H}$ at level $\alpha$ for the nuisance parameter value $\gamma$ when testing $H_{0}$ against $H_{1}$. \end{thm}
    
Based on this two-sided test, we can construct an asymptotically optimal $100(1-\alpha)$\% confidence set for $\theta$ by the test inversion:
\begin{equation*}
	CS=\left\{ \theta\in\Theta:Y_{i}\leq g(X_{i},\theta+\hat{h}_{n}^{-}(\theta)/n)\text{ for some }i\text{ and }Y_{i}\geq g(X_{i},\theta)\text{ for all }i\right\},
\end{equation*}
where $\hat{h}_{n}^{-}(\vartheta)$ is constructed at each $\vartheta$ from the plug-in estimator of $\lambda$ stated before Definition \ref{def:AUMPnui}, with $\theta_{0}=\vartheta$.

\section{Simulation\label{sec:sim}}

In this section, we investigate the finite-sample performance of the proposed tests through simulation.
We consider the symmetric first-price procurement auction with independent private costs discussed in Section \ref{sec:framework}. For $m=5$ bidders, we parameterize the conditional mean cost and the equilibrium bid boundary as $s(x,\theta)=\theta x-0.25$ and $g(x,\theta)=s(x,\theta)/(m-1)$, with $\theta_{0}=1.5$. The covariates $X_{i}$ are iid with support $\mathcal{X}=\{x_{1},x_{2},x_{3}\}=\{1,2,3\}$ and probabilities $(\pi_{1},\pi_{2},\pi_{3})=(0.3,0.4,0.3)$. Given a sample size $n$ and a local parameter $h$, we set $\theta=\theta_{0}+h/n$. Conditional on $X_{i}=x$, we draw the private cost $C_{i}$ from an exponential distribution with mean $s(x,\theta)$ and generate the observed equilibrium bid as $Y_{i} = C_{i}+g(X_{i},\theta)$, where $C_{i}\sim s(X_{i},\theta)\times\text{Exponential}(1)$.
Thus, conditional on $X_{i}=x$, the support of $Y_{i}$ is $[g(x,\theta),\infty)$. Thus, the density is a special case of the general model (\ref{eq:y|x,ga}) and is given by
\begin{equation*}
	f(y|x,\theta)
	=
	\frac{1}{s(x,\theta)}
	\exp\left(-\frac{y-g(x,\theta)}{s(x,\theta)}\right)
    \quad\text{and}\quad 
    g(x,\theta)=\frac{s(x,\theta)}{m-1}.
\end{equation*}
We set $\alpha=0.05$, consider $n\in\{20,40\}$, and use $10{,}000$ Monte Carlo replications.

In our proposed test, we use the test defined in (\ref{eq:phi-nui}) with $\hat{h}_{n}^{-}=\hat{\lambda}_{n}\log(1/\alpha)$ for the one-sided test $H_{0}^{-}:h\leq0$ against $H_{1}^{-}:h>0$ and the test defined in (\ref{eq:phi+nui}) with $\hat{h}_{n}^{+}=\hat{\lambda}_{n}\log(1/(1-\alpha))$ for the one-sided test $H_{0}^{+}:h\geq0$ against $H_{1}^{+}:h<0$, where $\hat{\lambda}_{n}$ is estimated by
\begin{equation*}
	\hat{\lambda}_{n} = \left(\frac{1}{n}\sum_{i=1}^{n}f(g(X_{i},\theta_{0})|X_{i},\theta_{0})\nabla_{\theta}g(X_{i},\theta_{0})\right)^{-1}.
\end{equation*}

In the CH test, we reject $H_{0}^{-}$ if $n(\hat{\theta}-\theta_{0})>q_{1-\alpha}(Z^{\theta_{0}})$, and reject $H_{0}^{+}$ if $n(\hat{\theta}-\theta_{0})<q_{\alpha}(Z^{\theta_{0}})$, where $\hat{\theta}$ denotes the maximum likelihood estimator. Since the lower-tail critical value $q_{\alpha}(Z^{\theta_{0}})$ leads to over-rejection, we also consider a modified rule that rejects $H_{0}^{+}$ if $n(\hat{\theta}-\theta_{0})<q_{0}(Z^{\theta_{0}})=0$.
Under $P_{\theta_{0}+h/n}$, the MLE satisfies $n(\hat{\theta}-(\theta_{0}+h/n))\rightsquigarrow Z^{\theta_{0}}$, where $q_{1-\alpha}(Z^{\theta_{0}})$ is the $(1-\alpha)$-th quantile of $Z^{\theta_{0}}$, and
\begin{equation*}
	Z^{\theta_{0}}=\frac{\mathrm{Exp}(1)}{E_{X}[\nabla_{\theta}g(X,\theta_{0})/s(X,\theta_{0})]}.
\end{equation*}
For feasible implementation, we replace the population expectation $E_{X}[\nabla_{\theta}g(X,\theta_{0})/s(X,\theta_{0})]$ with its sample analog, $\hat{\lambda}_{n}^{-1}$, to calculate the critical values in the CH test.

\begin{figure}[H]
	\centering

	\subfigure[$n=20$]{\centering\includegraphics[scale=0.25]{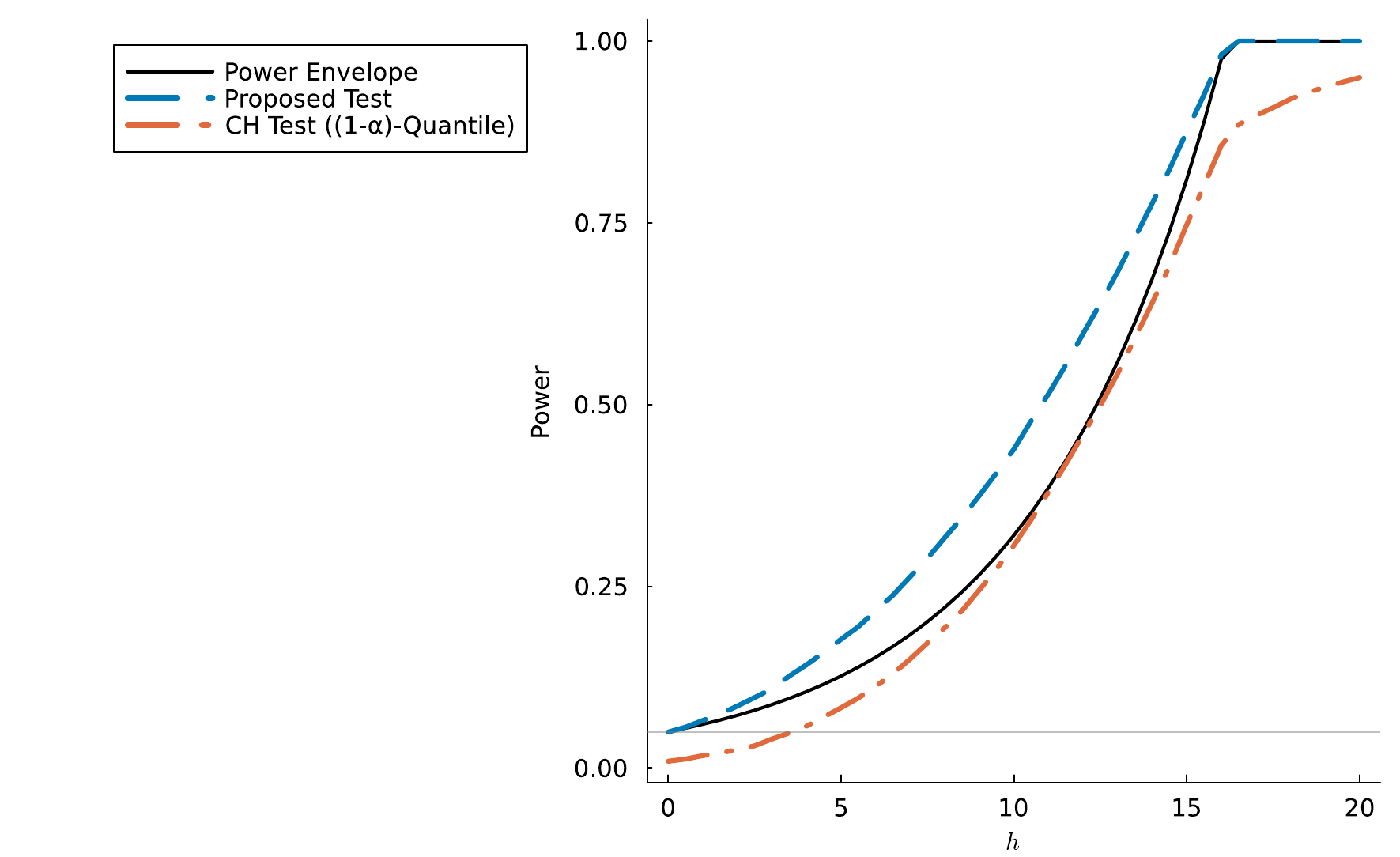}}\subfigure[$n=40$]{\centering\includegraphics[scale=0.25]{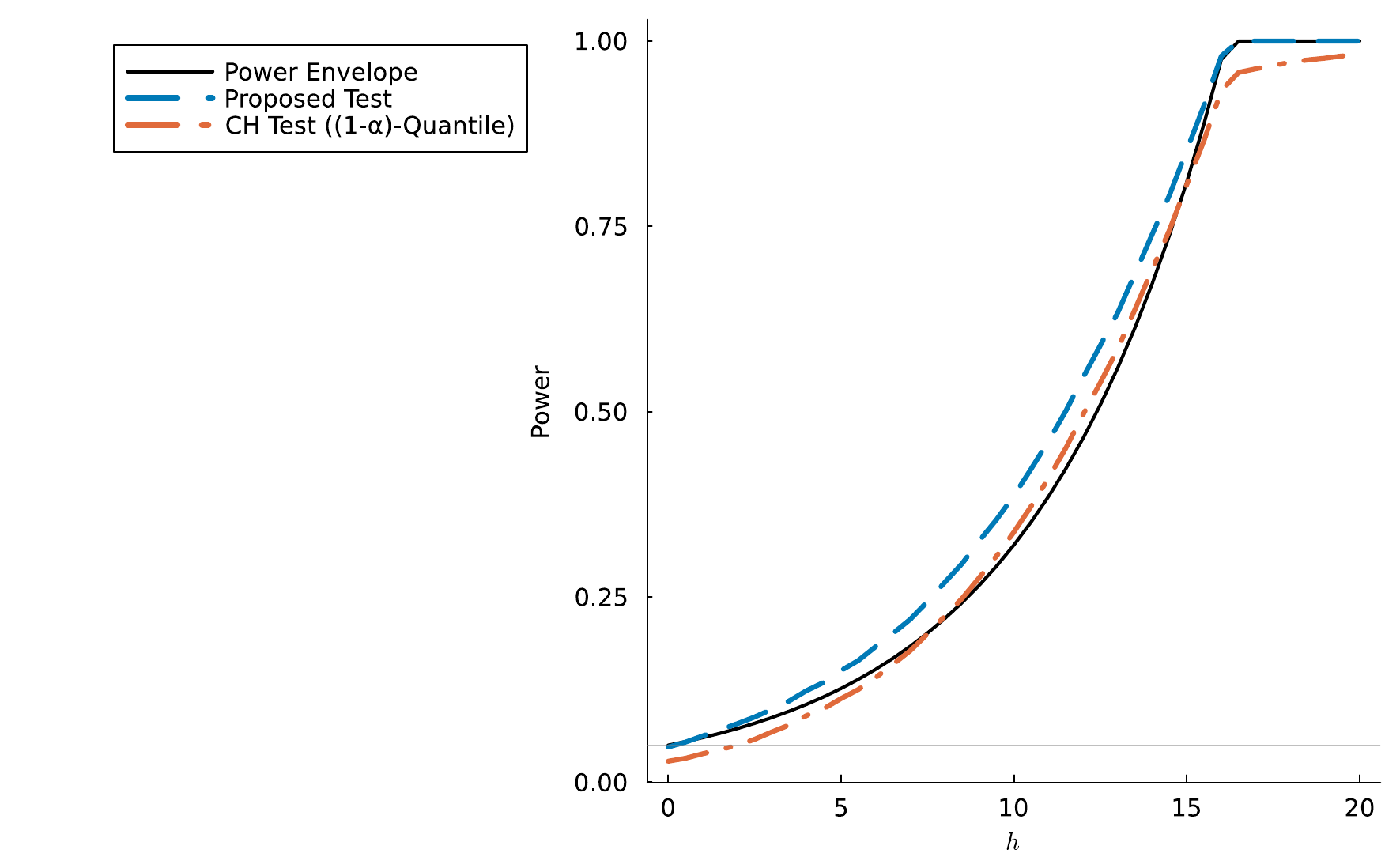}}
	\caption{Comparison of the proposed test with the CH test (the Wald test of \citet{chernozhukov2004likelihood}) at sample sizes $n\in\{20,40\}$ for $H_{0}^{-}:h\protect\leq0$ against $H_{1}^{-}:h>0$.}
	\label{fig:ch_cov_u}
\end{figure}

\begin{figure}[H]
	\centering

	\subfigure[$n=20$]{\centering\includegraphics[scale=0.25]{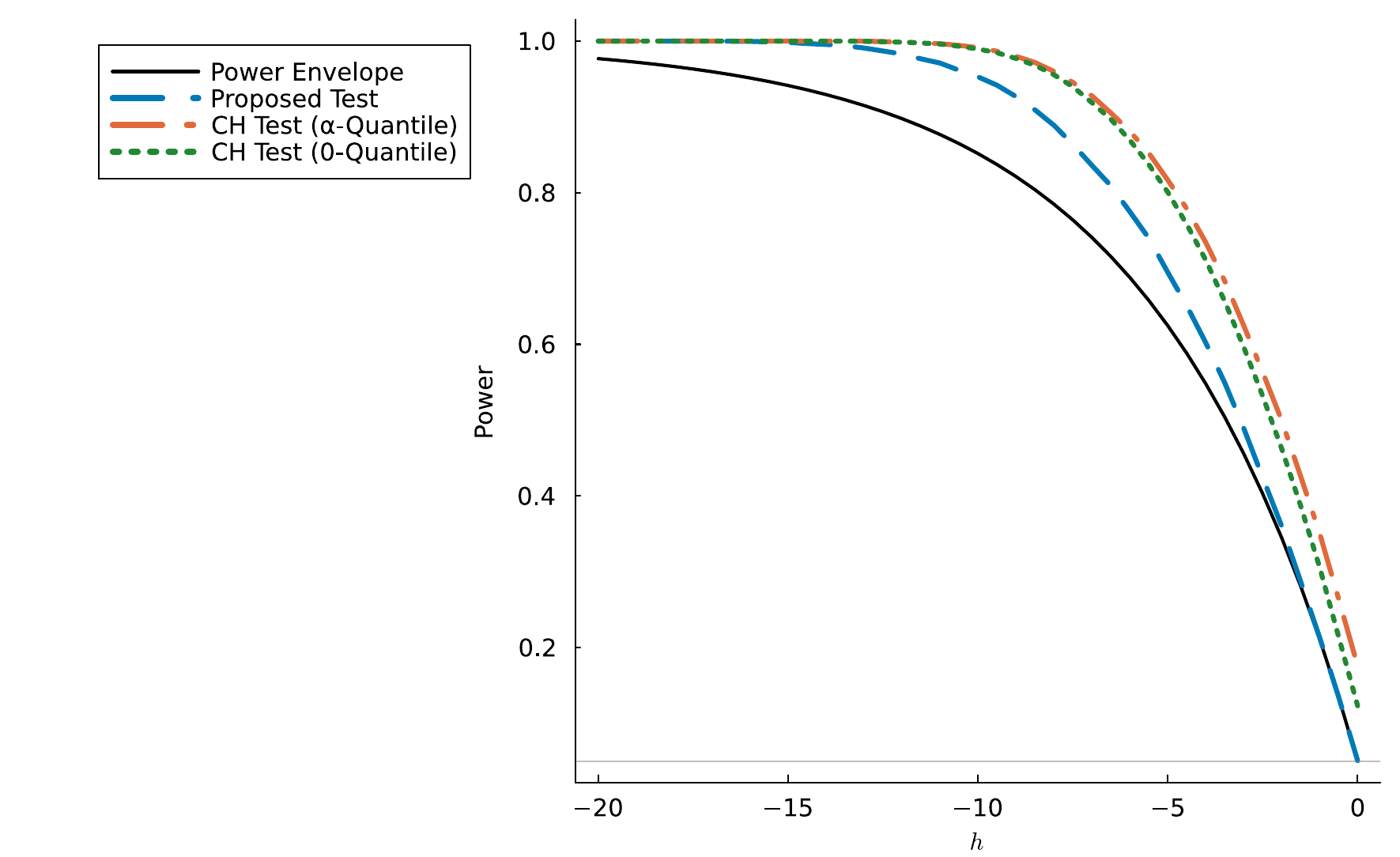}}\subfigure[$n=40$]{\centering\includegraphics[scale=0.25]{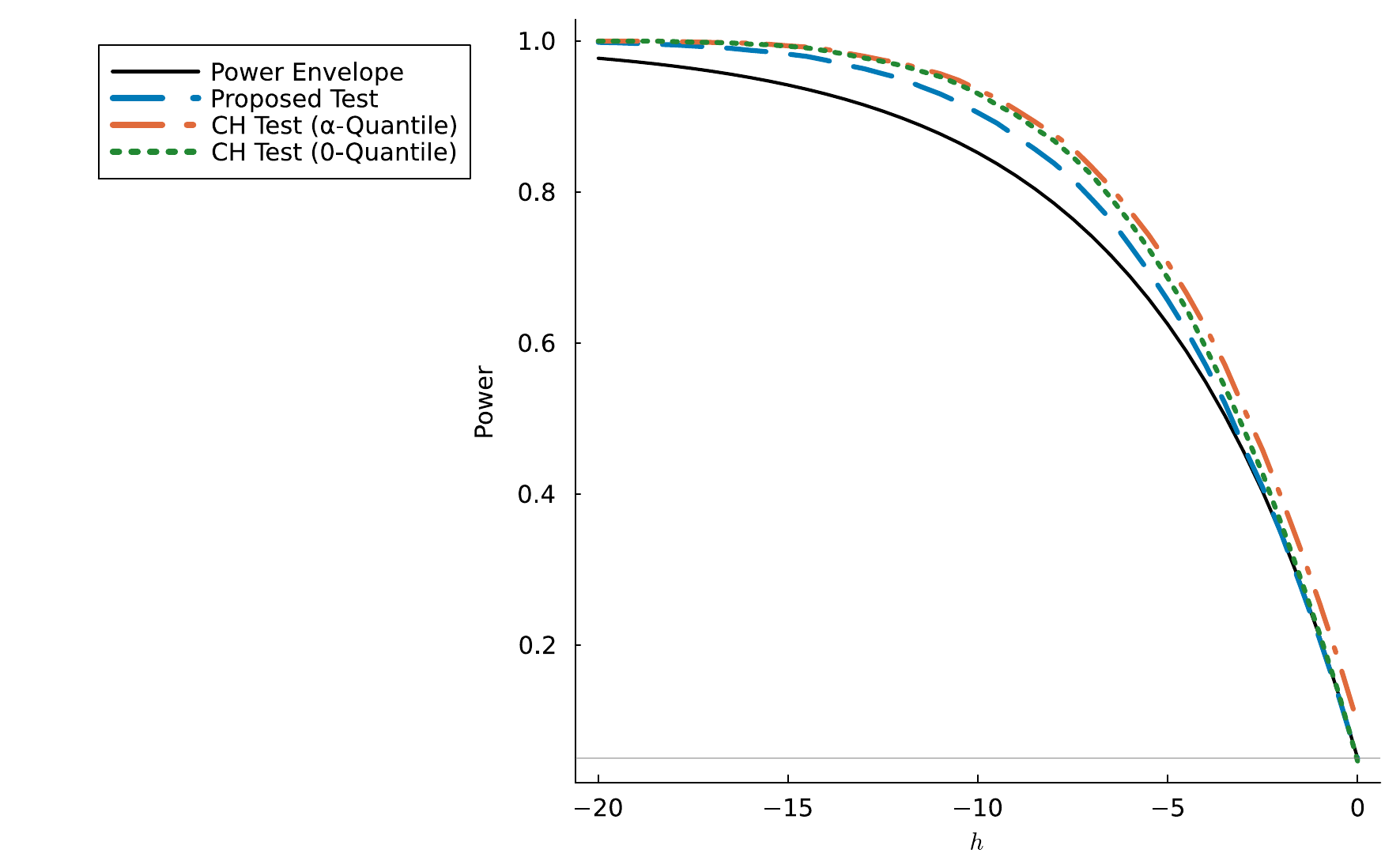}}
	\caption{Comparison of the proposed test with the CH test (the Wald test of \citet{chernozhukov2004likelihood}) at sample sizes $n\in\{20,40\}$ for $H_{0}^{+}:h\protect\geq0$ against $H_{1}^{+}:h<0$.}
	\label{fig:ch_cov_l}
\end{figure}

Figures \ref{fig:ch_cov_u} and \ref{fig:ch_cov_l} present the power curves of the proposed test and the CH test, together with the asymptotic power envelope for each one-sided testing problem.
The asymptotic power envelopes, derived in the proofs of Theorems \ref{thm:opt-nui} and \ref{thm:opt+nui} in the Appendix, represent the asymptotically optimal values.\footnote{The power envelope is given by $E_{h}[\phi^{-}(W_{h})]$ for the range $h\geq0$ and $E_{h}[\phi^{+}(W_{h})]$ for the range $h<0$ using the notation introduced in the Appendix.}
In Figure \ref{fig:ch_cov_u}, the proposed test rejects more frequently than the CH test over $h>0$, while the CH test is conservative near the null and large $h$. Both tests have valid size control. At $h=0$, the proposed test has sizes $4.95\%$ and $4.75\%$ for $n=20$ and $n=40$, respectively, compared with $0.94\%$ and $2.84\%$ for the CH test. 
In Figure \ref{fig:ch_cov_l}, the CH test exhibits severe over-rejection, whereas the proposed test has nearly valid size control. At $h=0$, the proposed test has sizes $5.07\%$ and $4.92\%$ for $n=20$ and $n=40$, respectively; the corresponding sizes are $17.57\%$ and $9.61\%$ for the CH test based on the $\alpha$-quantile and $12.31\%$ and $4.60\%$ for the CH test based on the zero quantile.

\appendix

\section{Mathematical appendix\label{app:math}}

\textbf{Notation:} Let $P_{h}\{\cdot\}$ and $E_{h}[\cdot]$ be the probability and the expectation under $f_{W}(w|h)$, respectively.

\subsection{Proof of Lemma \ref{lem:weak_conv}}

We omit the proof since it is a special case of the more general result in Lemma \ref{lem:weak_conv_nui}.

\subsection{Proof of Theorem \ref{thm:opt-}}

\subsubsection*{Step 1: Derive the UMP test in the limit of experiments}

Let $P_{h}^{W}$ denote the distribution of $W_{h}$, and write $\mathcal{P}^{W}=\{P_{h}^{W}:h\in\mathcal{H}\}$. For every $h_{1}<h_{2}$, $f_{W}(w|h_{1})=0$ implies $f_{W}(w|h_{2})=0$ for every $w$. On the support of $P_{h_{1}}^{W}$, the likelihood ratio is
\begin{equation*}
	\frac{\mathrm{d}P_{h_{2}}^{W}}{\mathrm{d}P_{h_{1}}^{W}}(w)
	=\frac{f_{W}(w|h_{2})}{f_{W}(w|h_{1})}
	=e^{(h_{2}-h_{1})/\lambda}\mathbb{I}\{w>h_{2}\}.
\end{equation*}
This likelihood ratio equals zero for $h_{1}<w\leq h_{2}$ and $e^{(h_{2}-h_{1})/\lambda}$ for $w>h_{2}$, and is therefore nondecreasing in $w$. Hence, $\mathcal{P}^{W}$ has monotone likelihood ratio in the statistic $S(W)=W$.

Let $\bar{h}^{-}=\lambda\log(1/\alpha)$. Observe that for every $\bar{h}$, we have
\begin{equation*}
	P_{h}\{W_{h}>\bar{h}\}=\frac{1}{\lambda}\int_{\max\{\bar{h},h\}}^{\infty}e^{-(w-h)/\lambda}\mathrm{d}w=\min\{e^{(h-\bar{h})/\lambda},1\},
\end{equation*}
which implies $P_{0}\{W_{0}>\bar{h}^{-}\}=\alpha$. Combining this with Lemma $\text{\ref{lem:MLRT}}$ implies that the test $\phi^{-}(W_{h}):=\mathbb{I}\{W_{h}>\bar{h}^{-}\}$ is UMP for testing $H_{0}:h\leq0$ against $H_{1}:h>0$ in the limit of experiments. Then the proposed test $\phi_{n}^{-}(Y^{n})$ in (\ref{eq:phi-}) emerges as its sample counterpart.

\subsubsection*{Step 2: Show size control of $\phi_{n}^{-}(Y^{n})$}

Pick an arbitrary $h\leq0$. Since the weak convergence in Lemma $\ref{lem:weak_conv}$ implies
\begin{equation*}
	\lim_{n}E_{\theta_{0}+\frac{h}{n}}[\phi_{n}^{-}(Y^{n})]=E_{h}[\phi^{-}(W_{h})],
\end{equation*}
we obviously have
\begin{align*}
	\lim_{n}E_{\theta_{0}+\frac{h}{n}}[\phi_{n}^{-}(Y^{n})] &=\alpha e^{h/\lambda}\leq\alpha,
\end{align*}
i.e., $\phi_{n}^{-}(Y^{n})$ achieves asymptotic size control.

\subsubsection*{Step 3: Show power optimality of $\phi_{n}^{-}(Y^{n})$}

Pick an arbitrary sequence of tests $\{\psi_{n}\}$ satisfying $\limsup_{n}E_{\theta_{0}}[\psi_{n}]\leq\alpha$, and pick an arbitrary $\bar{h}>0$. Consider a subsequence $\{\psi_{m}\}$ of $\{\psi_{n}\}$ such that
\begin{eqnarray*}
    \lim_{m}E_{\theta_{0}}[\psi_{m}] & \leq & \limsup_{n}E_{\theta_{0}}[\psi_{n}],\\
    \lim_{m}E_{\theta_{0}+\frac{\bar{h}}{m}}[\psi_{m}] & = & \limsup_{n}E_{\theta_{0}+\frac{\bar{h}}{n}}[\psi_{n}].
\end{eqnarray*}
By Lemma $\ref{lem:vdv}$, if a sequence of power functions of tests $\{E_{\theta_{0}+\frac{h}{m}}[\psi_{m}]\}$ converges pointwise to some function (denoted by $\pi_{\psi}(h)$) for $h\in\{0,\bar{h}\}$, then $\pi_{\psi}(h)$ is a power function for testing $H_{0}:h=0$ against $H_{1}:h=\bar{h}$ in one sample $W_{h}\sim f_{W}(w|h)$ for $h\in\{0,\bar{h}\}$. On the other hand, the Neyman-Pearson lemma (e.g., Theorem 3.2.1 (ii) of \citealp{lehmann2022testing}) implies that $\phi^{-}(W_{h})$ is the most powerful test for this one-sample testing problem. Therefore, we have $E_{\bar{h}}[\phi^{-}(W_{\bar{h}})]\geq\pi_{\psi}(\bar{h})$. Repeating the argument for every $\bar{h}$, the conclusion is obtained.

\subsection{Proof of Theorem \ref{thm:opt+}}

\subsubsection*{Step 1: Derive the UMP test in the limit of experiments}

Consider the testing problem $H_{0}:h=0$ against $H_{1}:h=\bar{h}$ for every fixed $\bar{h}<0$ in the limit of experiments. Define $a(k)=P_{0}\{f_{W}(W_{0}|\bar{h})>kf_{W}(W_{0}|0)\}$, which is nonincreasing. Pick an arbitrary $\alpha\in(0,1)$ and let $k_{0}$ be such that $a(k_{0})\leq\alpha\leq a(k_{0}-0)$. Then the Neyman-Pearson test is defined as
\begin{equation*}
	\phi^{\mathrm{NP}}(W_{h})=\left\{ \begin{array}{lll}
		1                                           & \text{if} & f_{W}(W_{h}|\bar{h})>k_{0}f_{W}(W_{h}|0) \\
		\frac{\alpha-a(k_{0})}{a(k_{0}-0)-a(k_{0})} & \text{if} & f_{W}(W_{h}|\bar{h})=k_{0}f_{W}(W_{h}|0) \\
		0                                           & \text{if} & f_{W}(W_{h}|\bar{h})<k_{0}f_{W}(W_{h}|0)
	\end{array}\right.,
\end{equation*}
with $E_{0}[\phi^{\mathrm{NP}}(W_{0})]=\alpha$. By the Neyman-Pearson lemma, this test is most powerful. We can compute $k_{0}$ and $a(k_{0})$ from $E_{0}[\phi^{\mathrm{NP}}(W_{0})]=\alpha$ and show that
\begin{eqnarray*}
    \phi^{\mathrm{NP}}(W_{h}) & = & \left\{ \begin{array}{lll}
        1      & \text{if} & f_{W}(W_{h}|\bar{h})>e^{\bar{h}/\lambda}f_{W}(W_{h}|0) \\
        \alpha & \text{if} & f_{W}(W_{h}|\bar{h})=e^{\bar{h}/\lambda}f_{W}(W_{h}|0) \\
        0      & \text{if} & f_{W}(W_{h}|\bar{h})<e^{\bar{h}/\lambda}f_{W}(W_{h}|0)
    \end{array}\right.=\left\{ \begin{array}{lll}
        1      & \text{if} & \mathbb{I}\{W_{h}>\bar{h}\}>\mathbb{I}\{W_{h}>0\} \\
        \alpha & \text{if} & \mathbb{I}\{W_{h}>\bar{h}\}=\mathbb{I}\{W_{h}>0\} \\
        0      & \text{if} & \mathbb{I}\{W_{h}>\bar{h}\}<\mathbb{I}\{W_{h}>0\}
    \end{array}\right.\\
    & = & \begin{cases}
        1      & \text{ if }\bar{h}<W_{h}\leq0 \\
        \alpha & \text{ if }W_{h}>0 \text{ or }W_{h}\leq\bar{h}
    \end{cases}
\end{eqnarray*}
is most powerful, where the second equality follows from the definition of $f_{W}$, and the third equality follows from the fact that $\mathbb{I}\{W_{h}>0\}=1$ implies $\mathbb{I}\{W_{h}>\bar{h}\}=1$ for every $\bar{h}<0$.

Define $\phi^{+}(W_{h})=(1-D_{h,0})+\alpha D_{h,0}$. Since $P_{0}\{W_{0}\leq\bar{h}\} = P_{\bar{h}}\{W_{\bar{h}}\leq\bar{h}\} = 0$, $\phi^{+}(W_{h})$ agrees with $\phi^{\mathrm{NP}}(W_{h})$ almost surely under both $P_{0}$ and $P_{\bar{h}}$.

Note that $\phi^{+}(\cdot)$ does not depend on the alternative $\bar{h}$. Thus, it is UMP for $H_{0}:h=0$ against $H_{1}:h<0$. Moreover, since $E_{h}[\phi^{+}(W_{h})]=\alpha$ for every $h>0$, it is also UMP for $H_{0}:h\geq0$ against $H_{1}:h<0$.

\subsubsection*{Step 2: Show size control of $\phi_{n}^{+}(Y^{n})$}

Let $\bar{h}^{+}=\lambda\log\left(\frac{1}{1-\alpha}\right)$. Observe that
\begin{eqnarray*}
    \lim_{n}E_{\theta_{0}+\frac{h}{n}}[\phi_{n}^{+}(Y^{n})] & = & \lim_{n}P_{\theta_{0}+\frac{h}{n}}\left\{ Y_{(1)}<\max\left\{ g(\theta_{0}),g\left(\theta_{0}+\frac{\bar{h}^{+}}{n}\right)\right\} \right\} \\
    & = & \lim_{n}P_{\theta_{0}+\frac{h}{n}}\left\{ Y_{(1)}<g\left(\theta_{0}+\frac{\bar{h}^{+}}{n}\right)\right\} =P_{h}\{W_{h}\leq\bar{h}^{+}\}\\
    & = & 1-\min\{(1-\alpha)e^{h/\lambda},1\}=\max\{1-(1-\alpha)e^{h/\lambda},0\},
\end{eqnarray*}
where the second equality follows because $g(\theta_{0})<g\left(\theta_{0}+\frac{\bar{h}^{+}}{n}\right)$ eventually, and the third equality follows from the weak convergence in Lemma $\ref{lem:weak_conv}$. Thus, we have $\lim_{n}E_{\theta_{0}+\frac{h}{n}}[\phi_{n}^{+}(Y^{n})]=\max\{1-(1-\alpha)e^{h/\lambda},0\}\leq\alpha$ for every $h\geq0$.

\subsubsection*{Step 3: Show power optimality of $\phi_{n}^{+}(Y^{n})$}

We can derive
\begin{equation*}
	E_{h}[\phi^{+}(W_{h})]=1-\min\{e^{h/\lambda},1\}+\alpha\min\{e^{h/\lambda},1\}=\max\{1-(1-\alpha)e^{h/\lambda},\alpha\},
\end{equation*}
and
\begin{equation*}
	\lim_{n}E_{\theta_{0}+\frac{h}{n}}[\phi_{n}^{+}(Y^{n})]=E_{h}[\phi^{+}(W_{h})],
\end{equation*}
for each $h<0$. Therefore, the AUMP property of $\phi_{n}^{+}(Y^{n})$ follows from the same argument used to prove Theorem \ref{thm:opt-}.

\subsection{Proof of Theorem \ref{thm:opt_ts}}

\subsubsection*{Step 1: Show size control of $\phi_{n}(Y^{n})$}

Let $\bar{h}^{-}=\lambda\log\left(\frac{1}{\alpha}\right)$. For $h=0$, we have
\begin{eqnarray*}
    \lim_{n}E_{\theta_{0}}[\phi_{n}(Y^{n})] & = & \lim_{n}P_{\theta_{0}}\left\{ Y_{(1)}>g\left(\theta_{0}+\frac{\bar{h}^{-}}{n}\right)\text{ or }Y_{(1)}<g(\theta_{0})\right\} \\
    & \leq & \lim_{n}P_{\theta_{0}}\left\{ Y_{(1)}>g\left(\theta_{0}+\frac{\bar{h}^{-}}{n}\right)\right\} +\lim_{n}P_{\theta_{0}}\{Y_{(1)}<g(\theta_{0})\}\\
    & = & P_{0}\{W_{0}>\bar{h}^{-}\}+P_{0}\{W_{0}\leq0\}=e^{-\bar{h}^{-}/\lambda}+0\\
    & = & \alpha,
\end{eqnarray*}
where the inequality follows from the union bound, and the second equality follows from the weak convergence in Lemma $\ref{lem:weak_conv}$. Thus, $\phi_{n}(Y^{n})$ controls asymptotic size.

\subsubsection*{Step 2: Show power optimality of $\phi_{n}(Y^{n})$}

For every $h\neq0$, note that
\begin{eqnarray*}
    &  & \underset{n}{\lim}P_{\theta_{0}+\frac{h}{n}}\left\{ Y_{(1)}>g\left(\theta_{0}+\frac{\bar{h}^{-}}{n}\right)\text{ or }Y_{(1)}<g(\theta_{0})\right\} \\
    & = & 1-\underset{n}{\lim}P_{\theta_{0}+\frac{h}{n}}\left\{ g(\theta_{0})\leq Y_{(1)}\leq g\left(\theta_{0}+\frac{\bar{h}^{-}}{n}\right)\right\} \\
    & = & 1-\underset{n}{\lim}P_{\theta_{0}+\frac{h}{n}}\left\{ Y_{(1)}\leq g\left(\theta_{0}+\frac{\bar{h}^{-}}{n}\right)\right\} +\underset{n}{\lim}P_{\theta_{0}+\frac{h}{n}}\{Y_{(1)}<g(\theta_{0})\}\\
    & = & P_{h}\{W_{h}>\bar{h}^{-}\}+P_{h}\{W_{h}\leq0\},
\end{eqnarray*}
where the first and second equalities follow from the set relationship, and the third equality follows from the weak convergence in Lemma $\ref{lem:weak_conv}$. For $h>0$, we have
\begin{equation*}
	\underset{n}{\lim}P_{\theta_{0}+\frac{h}{n}}\left\{ Y_{(1)}>g\left(\theta_{0}+\frac{\bar{h}^{-}}{n}\right)\text{ or }Y_{(1)}<g(\theta_{0})\right\} =P_{h}\{W_{h}>\bar{h}^{-}\}+0=E_{h}[\phi^{-}(W_{h})].
\end{equation*}
On the other hand, for $h<0$, we have
\begin{equation*}
	\underset{n}{\lim}P_{\theta_{0}+\frac{h}{n}}\left\{ Y_{(1)}>g\left(\theta_{0}+\frac{\bar{h}^{-}}{n}\right)\text{ or }Y_{(1)}<g(\theta_{0})\right\} =\alpha e^{h/\lambda}+(1-e^{h/\lambda})=E_{h}[\phi^{+}(W_{h})].
\end{equation*}
Thus, the AUMP property follows from the same argument used to prove Theorem \ref{thm:opt-}.

Hereafter, let $P_{h}\{\cdot\}$ and $E_{h}[\cdot]$ be the probability and the expectation under $f_{W}(w|h)\equiv\prod_{j=1}^{L}f_{W_{j}}(w_{j}|h,\gamma)$. Define $P_{h,j}\{\cdot\}$ and $E_{h,j}[\cdot]$ as the corresponding probability and expectation under $f_{W_{j}}(w_{j}|h,\gamma)$.

\subsection{Proof of Lemma \ref{lem:weak_conv_nui}}
The proof follows a similar structure to the proof of Theorem 2 in \citet{hirano2003asymptotic}. The main difference is that the local parameter $\hat{h}_{n}$ is a sample-dependent random variable rather than being fixed, which requires additional arguments to handle the randomness in the local parameter, especially for establishing the weak convergence to a Bernoulli distribution.

It suffices to show ($\text{\ref{eq:D_ga}}$) and ($\text{\ref{eq:exp_ga}}$) since ($\text{\ref{eq:Z_ga_hat}}$) follows from Slutsky's lemma. Recall that
\begin{equation*}
	\prod_{i=1}^{n}\mathbb{I}\{Y_{i}\geq g(X_{i},\theta_{0}+h/n)\}=1
	\quad\text{under }P_{\theta_{0}+h/n,\gamma}.
\end{equation*}
Also, ($\text{\ref{eq:Z_ga_hat}}$) implies ($\text{\ref{eq:Z_ga}}$) by setting $\hat{h}_{n}=\bar{h}$ and $\hat{\gamma}_{n}=\gamma$, which satisfy the conditions in Assumption \ref{asm:nui}. The extension to any finite set of local parameters is straightforward.
Fix $h\in\mathcal{H}$, $\bar{h}\in I$, and $\gamma\in\Gamma$. Define
\begin{eqnarray*}
    D_{n} & = & \prod_{i=1}^{n}\mathbb{I}\{Y_{i}\geq g(X_{i},\theta_{0}+\hat{h}_{n}/n)\},\\
    R_{n} & = & \sum_{i=1}^{n}\{\log f(Y_{i}|X_{i},\theta_{0}+\hat{h}_{n}/n,\hat{\gamma}_{n})-\log f(Y_{i}|X_{i},\theta_{0}+h/n,\hat{\gamma}_{n})\}.
\end{eqnarray*}

We first show ($\text{\ref{eq:exp_ga}}$) by proving $R_{n}\stackrel{p}{\rightarrow}(\bar{h}-h)/\lambda$. Expand $R_{n}$ as
\begin{eqnarray*}
    R_{n} & = & \frac{1}{n}\sum_{i=1}^{n}\nabla_{\theta}\log f(Y_{i}|X_{i},\theta_{0}+\tilde{h}_{n}/n,\hat{\gamma}_{n})(\hat{h}_{n}-h)\\
    & = & \frac{1}{n}\sum_{i=1}^{n}\left\{ \nabla_{\theta}\log f(Y_{i}|X_{i},\theta_{0},\gamma)+\frac{\tilde{h}_{n}}{n}\nabla_{\theta\theta}\log f(Y_{i}|X_{i},\tilde{\theta}_{n},\tilde{\gamma}_{n})\right.\\
    &  & \qquad\left.+\sqrt{n}(\hat{\gamma}_{n}-\gamma)^{\prime}\frac{1}{\sqrt{n}}\nabla_{\gamma\theta}\log f(Y_{i}|X_{i},\tilde{\theta}_{n},\tilde{\gamma}_{n})\right\} (\hat{h}_{n}-h),
\end{eqnarray*}
where the first equality follows from an expansion of $\sum_i \log f(Y_{i}\mid X_{i},\theta_{0}+\hat{h}_{n}/n,\hat{\gamma}_{n})$ with $\tilde{h}_{n}$ being between $\hat{h}_{n}$ and $h$, and the second equality follows from an expansion of $\sum_i \nabla_{\theta}\log f(Y_{i}|X_{i},\theta_{0}+\tilde{h}_{n}/n,\hat{\gamma}_{n})$ with $(\tilde{\theta}_{n},\tilde{\gamma}_{n})$ being between $(\theta_{0}+\tilde{h}_{n}/n,\hat{\gamma}_{n})$ and $(\theta_{0},\gamma)$.

By Assumption \ref{asm:nui} (ii), we have $\sqrt{n}(\hat{\gamma}_{n}-\gamma)=O_{P_{\theta_{0}+\frac{h}{n},\gamma}}(1)$ and $\hat{h}_{n}-\bar{h}=o_{P_{\theta_{0}+\frac{h}{n},\gamma}}(1)$. Thus, $|\tilde{\theta}_{n} - \theta_{0}| \leq |\tilde{h}_{n}|/n \leq (|\hat{h}_{n} - \bar{h}| + |h| + |\bar{h}|)/n = o_{P_{\theta_{0}+\frac{h}{n},\gamma}}(1)$.
Assumption $\text{\ref{asm:nui}}$ (ii) and (iii) imply that $\frac{1}{n^{2}}\sum_{i=1}^{n}\nabla_{\theta\theta}\log f(Y_{i}|X_{i},\tilde{\theta}_{n},\tilde{\gamma}_{n})=o_{P_{\theta_{0}+\frac{h}{n},\gamma}}(1)$ and $\frac{1}{n^{3/2}}\sum_{i=1}^{n}\nabla_{\gamma\theta}\log f(Y_{i}|X_{i},\tilde{\theta}_{n},\tilde{\gamma}_{n})=o_{P_{\theta_{0}+\frac{h}{n},\gamma}}(1)$. Combining these results, we have
\begin{align*}
    R_{n}&=\frac{1}{n}\sum_{i=1}^{n}\nabla_{\theta}\log f(Y_{i}|X_{i},\theta_{0},\gamma)(\hat{h}_{n}-h)+o_{P_{\theta_{0}+\frac{h}{n},\gamma}}(1)\\
    &\stackrel{p}{\rightarrow}E_{\theta_{0},\gamma}[\nabla_{\theta}\log f(Y|X,\beta_{0})](\bar{h}-h),
\end{align*}
where the convergence follows from a triangular-array WLLN and the continuous mapping theorem. Also note that
\begin{eqnarray*}
    E_{\theta_{0},\gamma}[\nabla_{\theta}\log f(Y|X,\beta_{0})] & = & E_{X}\left[\int_{g(X,\theta_{0})}^{\infty}\nabla_{\theta}f(y|X,\beta_{0})dy\right]\\
    & = & E_{X}\left[\nabla_{\theta}\int_{g(X,\theta_{0})}^{\infty}f(y|X,\beta_{0})dy+f(g(X,\theta_{0})|X,\beta_{0})\nabla_{\theta}g(X,\theta_{0})\right]\\
    & = & E_{X}[f(g(X,\theta_{0})|X,\beta_{0})\nabla_{\theta}g(X,\theta_{0})]=\lambda^{-1},
\end{eqnarray*}
where the first equality follows from the law of iterated expectations and ($\text{\ref{eq:y|x,ga}}$), the second equality follows from Assumption $\text{\ref{asm:nui}}$ (iii) and the Leibniz integral rule, and the third equality follows from $E_{X}\left[\int_{g(X,\theta)}^{\infty}f(y|X,\beta)dy\right]=1$. From the continuous mapping theorem, ($\text{\ref{eq:exp_ga}}$) holds true.

Next, we show that $D_{n}\rightsquigarrow D_{h,\bar{h}}$ under $\ensuremath{P_{\theta_{0}+\frac{h}{n},\gamma}}$. Let $D_{n}^{-}=\prod_{i=1}^{n}\mathbb{I}\{Y_{i}\geq g(X_{i},\theta_{0}+\bar{h}/n+n^{-5/4})\}$ and $D_{n}^{+}=\prod_{i=1}^{n}\mathbb{I}\{Y_{i}\geq g(X_{i},\theta_{0}+\bar{h}/n-n^{-5/4})\}$. 
Assumption \ref{asm:nui} (ii) implies that 
\begin{equation*}
	P_{\theta_{0}+\frac{h}{n},\gamma}\left(\bar{h}/n-n^{-5/4}\leq \hat{h}_{n}/n\leq \bar{h}/n+n^{-5/4}\right)\rightarrow 1.
\end{equation*}
Define $A_n := \{D_n^- \leq D_n \leq D_n^+\}$. Since $\nabla_{\theta}g(X_{i},\theta_{0})>0$ by Assumption \ref{asm:nui} (iii), we have
\begin{equation*}
	P_{\theta_{0}+\frac{h}{n},\gamma}(A_n)\rightarrow 1.
\end{equation*}
Using $A_n$, we can bound $P_{\theta_{0}+\frac{h}{n},\gamma}(D_n \leq z)$ from above and below using set relationships:
\begin{align*}
    P_{\theta_{0}+\frac{h}{n},\gamma}(D_n \leq z)
    &\leq P_{\theta_{0}+\frac{h}{n},\gamma}(\{D_n \leq z\} \cap A_n) + P_{\theta_{0}+\frac{h}{n},\gamma}(A_n^c) \\
    &\leq P_{\theta_{0}+\frac{h}{n},\gamma}(\{D_n^- \leq z\} \cap A_n) + P_{\theta_{0}+\frac{h}{n},\gamma}(A_n^c) \\
    &\leq P_{\theta_{0}+\frac{h}{n},\gamma}(D_n^- \leq z) + P_{\theta_{0}+\frac{h}{n},\gamma}(A_n^c),
\end{align*}
and
\begin{align*}
    P_{\theta_{0}+\frac{h}{n},\gamma}(D_n \leq z)
    &\geq P_{\theta_{0}+\frac{h}{n},\gamma}(\{D_n \leq z\} \cap A_n) \\
    &\geq P_{\theta_{0}+\frac{h}{n},\gamma}(\{D_n^+ \leq z\} \cap A_n) \\
    &\geq P_{\theta_{0}+\frac{h}{n},\gamma}(D_n^+ \leq z) - P_{\theta_{0}+\frac{h}{n},\gamma}(A_n^c).
\end{align*}
Therefore, it suffices to show $D_n^- \rightsquigarrow D_{h,\bar{h}}$ and $D_n^+ \rightsquigarrow D_{h,\bar{h}}$ under $\ensuremath{P_{\theta_{0}+\frac{h}{n},\gamma}}$ to show $D_n \rightsquigarrow D_{h,\bar{h}}$ under $\ensuremath{P_{\theta_{0}+\frac{h}{n},\gamma}}$.

Let $\bar{\theta}_{n}=\theta_{0}+\bar{h}/n+n^{-5/4}$. Under $\ensuremath{P_{\theta_{0}+\frac{h}{n},\gamma}}$, we have
\begin{equation*}
	D_{n}^{-}=\prod_{i=1}^{n}\mathbb{I}\{Y_{i}\geq g(X_{i},\bar{\theta}_{n})\}\rightsquigarrow\mathrm{Bernoulli}\left(\lim_{n\rightarrow\infty}E_{\theta_{0}+\frac{h}{n},\gamma}\left[\prod_{i=1}^{n}\mathbb{I}\{Y_{i}\geq g(X_{i},\bar{\theta}_{n})\}\right]\right).
\end{equation*}
We now compute $\lim_{n\rightarrow\infty}E_{\theta_{0}+\frac{h}{n},\gamma}\left[\prod_{i=1}^{n}\mathbb{I}\{Y_{i}\geq g(X_{i},\bar{\theta}_{n})\}\right]$. Note that, for each $j=1,\ldots,L$ and all sufficiently large $n$, there are deterministic $\grave{h}_{n,j}$ and $\tilde{\theta}_{n,j}$ such that $\theta_{0}+\grave{h}_{n,j}/n$ and $\tilde{\theta}_{n,j}$ are between $\theta_{0}+\bar{h}/n+n^{-5/4}$ and $\theta_{0}+h/n$. Applying the mean-value theorem separately for each $j$, we obtain
\begin{eqnarray*}
	&  & E_{\theta_{0}+\frac{h}{n},\gamma}[\mathbb{I}\{Y\geq g(X,\bar{\theta}_{n})\}]=1-E_{X}[F_{Y\mid X,\theta_{0}+\frac{h}{n},\gamma}(g(X,\bar{\theta}_{n}))]\\
	& = & 1-E_{X}\left[\int_{-\infty}^{g(X,\bar{\theta}_{n})}f(y|X,\theta_{0}+h/n,\gamma)\mathbb{I}\{y-g(X,\theta_{0}+h/n)\geq0\}dy\right]\\
	& = & 1-E_{X}\left[\mathbb{I}\{g(X,\bar{\theta}_{n})-g(X,\theta_{0}+h/n)\geq0\}\int_{g(X,\theta_{0}+h/n)}^{g(X,\bar{\theta}_{n})}f(y|X,\theta_{0}+h/n,\gamma)dy\right]\\
	& = & 1-\sum_{j=1}^{L}\Pr\{X=a_{j}\}\left[\begin{array}{c}
			\mathbb{I}\{g(a_{j},\bar{\theta}_{n})-g(a_{j},\theta_{0}+h/n)\geq0\}\\
			\times f(g(a_{j},\theta_{0}+\grave{h}_{n,j}/n)|a_{j},\theta_{0}+h/n,\gamma)\\
			\times\{g(a_{j},\bar{\theta}_{n})-g(a_{j},\theta_{0}+h/n)\}
		\end{array}\right]\\
	& = & 1-\frac{1}{n}\sum_{j=1}^{L}\Pr\{X=a_{j}\}\left[\begin{array}{c}
			\mathbb{I}\{\nabla_{\theta}g(a_{j},\tilde{\theta}_{n,j})(\bar{h}+n^{-1/4}-h)\geq0\}\\
			\times f(g(a_{j},\theta_{0}+\grave{h}_{n,j}/n)|a_{j},\theta_{0}+h/n,\gamma)\\
			\times\nabla_{\theta}g(a_{j},\tilde{\theta}_{n,j})(\bar{h}+n^{-1/4}-h)
			\end{array}\right],
\end{eqnarray*}
where the second equality follows from the model ($\text{\ref{eq:y|x,ga}}$), the fourth equality follows from the integral mean-value theorem and the local strict monotonicity of $g(a_{j},\cdot)$ applied separately for each $j$, and the fifth equality follows from a mean-value expansion of $g(a_{j},\bar{\theta}_{n})$ around $\theta_{0}+h/n$ for each $j$. Thus, recalling that the intermediate points are non-random for every $j$, we have
\begin{eqnarray*}
	&  & E_{\theta_{0}+\frac{h}{n},\gamma}\left[\prod_{i=1}^{n}\mathbb{I}\{Y_{i}\geq g(X_{i},\bar{\theta}_{n})\}\right]=[P_{\theta_{0}+\frac{h}{n},\gamma}\{Y\geq g(X,\bar{\theta}_{n})\}]^{n}\\
	& = & \left[1-\frac{1}{n}\sum_{j=1}^{L}\Pr\{X=a_{j}\}\left[\begin{array}{c}
				\mathbb{I}\{\nabla_{\theta}g(a_{j},\tilde{\theta}_{n,j})(\bar{h}+n^{-1/4}-h)\geq0\}\\
				\times f(g(a_{j},\theta_{0}+\grave{h}_{n,j}/n)|a_{j},\theta_{0}+h/n,\gamma)\\
				\times\nabla_{\theta}g(a_{j},\tilde{\theta}_{n,j})(\bar{h}+n^{-1/4}-h)
			\end{array}\right]\right]^{n}\\
	& \rightarrow & \exp\left(\lim_{n\rightarrow\infty}-\sum_{j=1}^{L}\Pr\{X=a_{j}\}\left[\begin{array}{c}
			\mathbb{I}\{\nabla_{\theta}g(a_{j},\tilde{\theta}_{n,j})(\bar{h}+n^{-1/4}-h)\geq0\}\\
			\times f(g(a_{j},\theta_{0}+\grave{h}_{n,j}/n)|a_{j},\theta_{0}+h/n,\gamma)\\
			\times\nabla_{\theta}g(a_{j},\tilde{\theta}_{n,j})(\bar{h}+n^{-1/4}-h)
		\end{array}\right]\right)\\
	& = & \exp\left(-\sum_{j=1}^{L}\Pr\{X=a_{j}\}\lim_{n\rightarrow\infty}\left[\begin{array}{c}
			\mathbb{I}\{\nabla_{\theta}g(a_{j},\tilde{\theta}_{n,j})(\bar{h}+n^{-1/4}-h)\geq0\}\\
			\times f(g(a_{j},\theta_{0}+\grave{h}_{n,j}/n)|a_{j},\theta_{0}+h/n,\gamma)\\
			\times\nabla_{\theta}g(a_{j},\tilde{\theta}_{n,j})(\bar{h}+n^{-1/4}-h)
		\end{array}\right]\right)\\
	& = & \exp\left(-\sum_{j=1}^{L}\Pr\{X=a_{j}\}\left[\begin{array}{c}
			\lim_{n\rightarrow\infty}\mathbb{I}\{\nabla_{\theta}g(a_{j},\tilde{\theta}_{n,j})(\bar{h}+n^{-1/4}-h)\geq0\}\\
			\times f(g(a_{j},\theta_{0})|a_{j},\theta_{0},\gamma)\\
			\times\nabla_{\theta}g(a_{j},\theta_{0})(\bar{h}-h)
		\end{array}\right]\right)\\
	& = & \exp\left(-E_{X}\left[\begin{array}{c}
			\mathbb{I}\{\bar{h}-h\geq0\}\\
			\times f(g(X,\theta_{0})|X,\theta_{0},\gamma)\nabla_{\theta}g(X,\theta_{0})(\bar{h}-h)
		\end{array}\right]\right),
\end{eqnarray*}
where the first equality follows from the iid assumption, the second equality follows from the calculation above, the convergence follows from the standard exponential limit, the third equality follows from the finiteness of $\mathcal{X}$, the fourth equality follows from the continuity of $f(y|x,\theta,\gamma)$ in $y$ and $\theta$ and the continuity of $g(x,\theta)$ and $\nabla_{\theta}g(x,\theta)$ in $\theta$, and the last equality follows from the indicator limit established below. The assumption $\nabla_{\theta}g(a_{j},\theta_{0})>0$ guarantees
\begin{eqnarray*}
	&  & \lim_{n}\mathbb{I}\{\nabla_{\theta}g(a_{j},\tilde{\theta}_{n,j})(\bar{h}+n^{-1/4}-h)\geq0\}\\
	& = & \lim_{n}\mathbb{I}\{\bar{h}+n^{-1/4}-h\geq0\}=\ensuremath{\mathbb{I}\left[\bigcap_{n=1}^{\infty}\{\bar{h}-h\geq-n^{-1/4}\}\right]}=\mathbb{I}\{\bar{h}-h\geq0\},
\end{eqnarray*}
for each $j=1,\ldots,L$. Therefore, we obtain
\begin{eqnarray*}
	D_{n}^{-} & \rightsquigarrow & \mathrm{Bernoulli}(\exp(-E_{X}[\mathbb{I}\{\bar{h}-h\geq0\}f(g(X,\theta_{0})|X,\theta_{0},\gamma)\nabla_{\theta}g(X,\theta_{0})(\bar{h}-h)]))\\
	& = & D_{h,\bar{h}}.
\end{eqnarray*}
Similarly, we obtain $D_{n}^{+}\rightsquigarrow D_{h,\bar{h}}$. Combining these results yields ($\text{\ref{eq:D_ga}}$). Therefore, the conclusion is obtained.

\subsection{Proof of Theorem \ref{thm:opt-nui}}

\subsubsection*{Step 1: Derive the UMP test in the limit of experiments}

The proof of this part is similar to the one for Theorem $\text{\ref{thm:opt-}}$ after replacing $Z(h,\bar{h})$, $D_{h,\bar{h}}$, and $\lambda$ in Section $\text{\ref{sec:bench}}$ with the ones in Section $\text{\ref{sec:gen}}$ and redefining $f_{W}(w|h):=\prod_{j=1}^{L}f_{W_{j}}(w_{j}|h,\gamma)$. The only differences are as follows. For each $h_{1}<h_{2}$, the support of $f_{W}(\cdot|h_{2})$ is contained in the support of $f_{W}(\cdot|h_{1})$. On the support of $f_{W}(\cdot|h_{1})$,
\begin{align*}
	\frac{f_{W}(w|h_{2})}{f_{W}(w|h_{1})}
	&=\exp\left((h_{2}-h_{1})\sum_{j=1}^{L}\frac{G_{j}}{\lambda_{j}}\right)
	\mathbb{I}\left\{\min_{j\in\{1,\ldots,L\}}\frac{w_{j}}{G_{j}}>h_{2}\right\}\\
	&=\exp\left(\frac{h_{2}-h_{1}}{\lambda}\right)
	\mathbb{I}\left\{\min_{j\in\{1,\ldots,L\}}\frac{w_{j}}{G_{j}}>h_{2}\right\}.
\end{align*}
Thus, the likelihood ratio is nondecreasing in $\min_{j\in\{1,\ldots,L\}}w_{j}/G_{j}$, and the family has monotone likelihood ratio in this statistic. Observe that
\begin{equation*}
	P_{h,j}\{W_{h,j}>G_{j}\bar{h}\}=\frac{1}{\lambda_{j}}\int_{\max\{G_{j}\bar{h},G_{j}h\}}^{\infty}e^{-(w_{j}-G_{j}h)/\lambda_{j}}\mathrm{d}w_{j}=\exp\left(\min\left\{ \frac{G_{j}(h-\bar{h})}{\lambda_{j}},0\right\} \right).
\end{equation*}

Hence,
\begin{eqnarray*}
    &  & P_{h}\{\min_{j\in\{1,\ldots,L\}}W_{h,j}/G_{j}>\bar{h}\}\\
    & = & \prod_{j=1}^{L}P_{h,j}\{W_{h,j}>G_{j}\bar{h}\}=\prod_{j=1}^{L}\exp\left(\min\left\{ \frac{G_{j}(h-\bar{h})}{\lambda_{j}},0\right\} \right)\\
    & = & \exp\left(\min\left\{ \sum_{j=1}^{L}\frac{G_{j}(h-\bar{h})}{\lambda_{j}},0\right\} \right)=\min\left\{ \exp\left(\frac{h-\bar{h}}{\lambda}\right),1\right\} ,
\end{eqnarray*}
where the first equality holds since $(W_{h,1},\ldots,W_{h,L})$ are mutually independent, the second equality follows from the above observation, the third equality from $G_{j}>0$ and $\lambda_{j}>0$ for every $j$, and the last equality from $\lambda=(\sum_{j=1}^{L}G_{j}/\lambda_{j})^{-1}$. Thus, $P_{0}\{\min_{j\in\{1,\ldots,L\}}W_{0,j}/G_{j}>\bar{h}^{-}\}=\alpha$ since $\bar{h}^{-}=\lambda\log(1/\alpha)$. Combining this with Lemma $\text{\ref{lem:MLRT}}$ implies that the test
\begin{equation*}
	\phi^{-}(W_{h}):=\mathbb{I}\left\{ \min_{j\in\{1,\ldots,L\}}W_{h,j}/G_{j}>\bar{h}^{-}\right\} =D_{h,\bar{h}^{-}}
\end{equation*}
is UMP for $H_{0}:h\leq0$ versus $H_{1}:h>0$ in the limit experiment. Then the proposed test $\phi_{n}^{-}(\hat{h}_{n}^{-},Y^{n},X^{n})$ in (\ref{eq:phi-nui}) emerges as its sample counterpart.

\subsubsection*{Step 2: Show size control of $\phi_{n}^{-}(\hat{h}_{n}^{-},Y^{n},X^{n})$}
Since the weak convergence in Lemma $\ref{lem:weak_conv_nui}$ implies\footnote{Under Assumption \ref{asm:nui}, both $\prod_{i=1}^{n}\mathbb{I}\{Y_{i}\geq g(X_{i},\theta_{0}+\hat{h}_{n}/n)\}$ and $\prod_{i=1}^{n}\mathbb{I}\{Y_{i}> g(X_{i},\theta_{0}+\hat{h}_{n}/n)\}$ have the same weak convergence property as in Lemma \ref{lem:weak_conv_nui}. Thus, the choice of strict or non-strict inequality does not affect the asymptotic results.}
\begin{equation*}
	\lim_{n}E_{\theta_{0}+\frac{h}{n},\gamma}[\phi_{n}^{-}(\hat{h}_{n}^{-},Y^{n},X^{n})]=E_{h}[\phi^{-}(W_{h})],
\end{equation*}
we obtain the asymptotic size control
\begin{align*}
    \lim_{n}E_{\theta_{0}+\frac{h}{n},\gamma}[\phi_{n}^{-}(\hat{h}_{n}^{-},Y^{n},X^{n})] & =\alpha e^{h/\lambda}\leq\alpha,
\end{align*}
for every $h\leq0$.

\subsubsection*{Step 3: Show power optimality of $\phi_{n}^{-}(\hat{h}_{n}^{-},Y^{n},X^{n})$}

This step follows from the same argument used to prove Theorem \ref{thm:opt-}, and the asymptotic power envelope is $E_{h}[\phi^{-}(W_{h})]=\min\{\alpha e^{h/\lambda},1\}$ for $h>0$. Therefore, the conclusion is obtained.

\subsection{Proof of Theorem \ref{thm:opt+nui}}

\subsubsection*{Step 1: Derive the UMP test in the limit of experiments}

The proof of this part is similar to the one for Theorem $\text{\ref{thm:opt+}}$ after replacing $D_{h,\bar{h}}$ and $\lambda$ in Section $\text{\ref{sec:bench}}$ with the ones in Section $\text{\ref{sec:gen}}$ and redefining $f_{W}(w|h):=\prod_{j=1}^{L}f_{W_{j}}(w_{j}|h,\gamma)$. The only differences are as follows. For every fixed $\bar{h}<0$,
\begin{eqnarray*}
    \phi^{+}(W_{h}) & = & \left\{ \begin{array}{lll}
        1      & \text{if} & f_{W}(W_{h}|\bar{h})>e^{\bar{h}/\lambda}f_{W}(W_{h}|0) \\
        \alpha & \text{if} & f_{W}(W_{h}|\bar{h})=e^{\bar{h}/\lambda}f_{W}(W_{h}|0) \\
        0      & \text{if} & f_{W}(W_{h}|\bar{h})<e^{\bar{h}/\lambda}f_{W}(W_{h}|0)
    \end{array}\right.\\
    & = & \left\{ \begin{array}{lll}
        1      & \text{if} & \mathbb{I}\{\min_{j\in\{1,\ldots,L\}}W_{h,j}/G_{j}>\bar{h}\}>\mathbb{I}\{\min_{j\in\{1,\ldots,L\}}W_{h,j}/G_{j}>0\} \\
        \alpha & \text{if} & \mathbb{I}\{\min_{j\in\{1,\ldots,L\}}W_{h,j}/G_{j}>\bar{h}\}=\mathbb{I}\{\min_{j\in\{1,\ldots,L\}}W_{h,j}/G_{j}>0\} \\
        0      & \text{if} & \mathbb{I}\{\min_{j\in\{1,\ldots,L\}}W_{h,j}/G_{j}>\bar{h}\}<\mathbb{I}\{\min_{j\in\{1,\ldots,L\}}W_{h,j}/G_{j}>0\}
    \end{array}\right.\\
    & = & \begin{cases}
        1      & \text{ if }\min_{j\in\{1,\ldots,L\}}W_{h,j}/G_{j}\leq0 \\
        \alpha & \text{ if }\min_{j\in\{1,\ldots,L\}}W_{h,j}/G_{j}>0
    \end{cases}=(1-D_{h,0})+\alpha D_{h,0}
\end{eqnarray*}
is most powerful. The second equality follows from simple transformations. The third equality follows from the fact that $\mathbb{I}\{\min_{j\in\{1,\ldots,L\}}W_{h,j}/G_{j}>0\}=1$ implies $\mathbb{I}\{\min_{j\in\{1,\ldots,L\}}W_{h,j}/G_{j}>\bar{h}\}=1$ for every $\bar{h}<0$ and holds almost surely under both $P_{0}$ and $P_{\bar{h}}$. The final equality follows from the definition of $D_{h,0}$.

\subsubsection*{Step 2: Show size control of $\phi_{n}^{+}(\hat{h}_{n}^{+},Y^{n},X^{n})$}

Observe that
\begin{align*}
	\lim_{n}E_{\theta_{0}+\frac{h}{n},\gamma}[\phi_{n}^{+}(\hat{h}_{n}^{+},Y^{n},X^{n})]
	&=\lim_{n}P_{\theta_{0}+\frac{h}{n},\gamma}
	\left\{\begin{array}{c}
		Y_{i}<\max\{g(X_{i},\theta_{0}),g(X_{i},\theta_{0}+\hat{h}_{n}^{+}/n)\}\\
		\text{for some }i
	\end{array}\right\}\\
	&=\lim_{n}P_{\theta_{0}+\frac{h}{n},\gamma}
	\{Y_{i}<g(X_{i},\theta_{0}+\hat{h}_{n}^{+}/n)\text{ for some }i\}\\
	&=1-\lim_{n}P_{\theta_{0}+\frac{h}{n},\gamma}
	\{Y_{i}\geq g(X_{i},\theta_{0}+\hat{h}_{n}^{+}/n)\text{ for all }i\}\\
	&=1-P_{h}\left\{\min_{j\in\{1,\ldots,L\}}W_{h,j}/G_{j}>\bar{h}^{+}\right\}\\
	&=\max\left\{1-\exp\left(\frac{h-\bar{h}^{+}}{\lambda}\right),0\right\}.
\end{align*}
where the second equality follows because $g(X_{i},\theta_{0})<g(X_{i},\theta_{0}+\hat{h}_{n}^{+}/n)$ with probability approaching one, and the fourth equality follows from the weak convergence in Lemma $\ref{lem:weak_conv_nui}$. Since $\bar{h}^{+}=\lambda\log(1/(1-\alpha))$, we obtain the asymptotic size control
\begin{align*}
	\lim_{n}E_{\theta_{0}+\frac{h}{n},\gamma}[\phi_{n}^{+}(\hat{h}_{n}^{+},Y^{n},X^{n})] &=\max\{1-(1-\alpha)e^{h/\lambda},0\}\leq\alpha,
\end{align*}
for every $h\geq0$.

\subsubsection*{Step 3: Show power optimality of $\phi_{n}^{+}(\hat{h}_{n}^{+},Y^{n},X^{n})$}

This step follows from the same argument used to prove Theorem \ref{thm:opt+}, and the asymptotic power envelope is $E_{h}[\phi^{+}(W_{h})]=\max\{1-(1-\alpha)e^{h/\lambda},0\}$ for $h<0$. Therefore, the conclusion is obtained.

\subsection{Proof of Theorem \ref{thm:opt_ts-nui}}

The proof of this part is similar to the one for Theorem $\text{\ref{thm:opt_ts}}$ after replacing $D_{h,\bar{h}}$ and $\lambda$ in Section $\text{\ref{sec:bench}}$ with the ones in Section $\text{\ref{sec:gen}}$.

\subsubsection*{Step 1: Show size control of $\phi_{n}(\hat{h}_{n}^{-},Y^{n},X^{n})$}

For $h=0$, we have
\begin{eqnarray*}
    &  & \limsup_{n}E_{\theta_{0},\gamma}[\phi_{n}(\hat{h}_{n}^{-},Y^{n},X^{n})]\\
    & = & \limsup_{n}P_{\theta_{0},\gamma}\left\{ Y_{i}>g(X_{i},\theta_{0}+\hat{h}_{n}^{-}/n)\text{ for all }i\text{ or }Y_{i}<g(X_{i},\theta_{0})\text{ for some }i\right\} \\
    & \leq & \lim_{n}P_{\theta_{0},\gamma}\left\{ Y_{i}>g(X_{i},\theta_{0}+\hat{h}_{n}^{-}/n)\text{ for all }i\right\} +\lim_{n}P_{\theta_{0},\gamma}\left\{ Y_{i}<g(X_{i},\theta_{0})\text{ for some }i\right\} \\
    & = & \exp\left(-\frac{\bar{h}^{-}}{\lambda}\right)+0=\alpha,
\end{eqnarray*}
where the inequality follows from the union bound, and the second equality follows from the weak convergence in Lemma \ref{lem:weak_conv_nui}. Thus, $\phi_{n}(\hat{h}_{n}^{-},Y^{n},X^{n})$ controls asymptotic size.

\subsubsection*{Step 2: Show power optimality of $\phi_{n}(\hat{h}_{n}^{-},Y^{n},X^{n})$}

For every $h\neq0$, note that
\begin{align*}
	&\underset{n}{\lim}P_{\theta_{0}+\frac{h}{n},\gamma}
	\left\{\begin{array}{c}
		Y_{i}>g(X_{i},\theta_{0}+\hat{h}_{n}^{-}/n)\text{ for all }i\\
		\text{or }Y_{i}<g(X_{i},\theta_{0})\text{ for some }i
	\end{array}\right\}\\
	&=1-\underset{n}{\lim}P_{\theta_{0}+\frac{h}{n},\gamma}
	\left\{\begin{array}{c}
		Y_{i}\leq g(X_{i},\theta_{0}+\hat{h}_{n}^{-}/n)\text{ for some }i\\
		\text{and }Y_{i}\geq g(X_{i},\theta_{0})\text{ for all }i
	\end{array}\right\}\\
	&=1-\left[
	\underset{n}{\lim}P_{\theta_{0}+\frac{h}{n},\gamma}
	\{Y_{i}\leq g(X_{i},\theta_{0}+\hat{h}_{n}^{-}/n)\text{ for some }i\}\right.\\
	&\qquad\left.
	-\underset{n}{\lim}P_{\theta_{0}+\frac{h}{n},\gamma}
	\{Y_{i}<g(X_{i},\theta_{0})\text{ for some }i\}
	\right]\\
	&=\underset{n}{\lim}P_{\theta_{0}+\frac{h}{n},\gamma}
	\{Y_{i}>g(X_{i},\theta_{0}+\hat{h}_{n}^{-}/n)\text{ for all }i\}\\
	&\qquad+\underset{n}{\lim}P_{\theta_{0}+\frac{h}{n},\gamma}
	\{Y_{i}<g(X_{i},\theta_{0})\text{ for some }i\},
\end{align*}
where the second equality follows from $g(X_{i},\theta_{0})\leq g(X_{i},\theta_{0}+\hat{h}_{n}^{-}/n)$ eventually, which implies the event $\{Y_{i}\leq g(X_{i},\theta_{0}+\hat{h}_{n}^{-}/n)\text{ for some }i\}$ includes $\{Y_{i}<g(X_{i},\theta_{0})\text{ for some }i\}$ eventually. For $h>0$, we have
\begin{eqnarray*}
    &  & \underset{n}{\lim}P_{\theta_{0}+\frac{h}{n},\gamma}\left\{ Y_{i}>g(X_{i},\theta_{0}+\hat{h}_{n}^{-}/n)\text{ for all }i\text{ or }Y_{i}<g(X_{i},\theta_{0})\text{ for some }i\right\} \\
    & = & \min\left\{ \exp\left(\frac{h-\bar{h}^{-}}{\lambda}\right),1\right\} +0=E_{h}[\phi^{-}(W_{h})].
\end{eqnarray*}
On the other hand, for $h<0$, we have
\begin{eqnarray*}
    &  & \underset{n}{\lim}P_{\theta_{0}+\frac{h}{n},\gamma}\left\{ Y_{i}>g(X_{i},\theta_{0}+\hat{h}_{n}^{-}/n)\text{ for all }i\text{ or }Y_{i}<g(X_{i},\theta_{0})\text{ for some }i\right\} \\
    & = & 1-(1-\alpha)e^{h/\lambda}=E_{h}[\phi^{+}(W_{h})].
\end{eqnarray*}
Thus, the AUMP property follows from the same argument used to prove Theorem \ref{thm:opt-}.

\bibliographystyle{apalike}
\bibliography{ref_list}

\end{document}